\documentclass{amsart}
\usepackage{amssymb}
\usepackage[hyperref]{degt}

\let\4\relax

\input NL.def

\def\prefix{6}
\mathchardef\minmin=260
\mathchardef\minmax=284
\let\mincount\minmin
\def\minusform#1#2{\frac{\bar{#1}}{#2}}

\def\sizemax{285}
\def\sizesub{261}

\def\itemrefform#1{{\rm(#1)}}

\def\2{\color{red}}

\makeatletter
\def\orbitlist{%
 \def\o@@{\ifnum\z@=\count\z@\single\else
   \ifnum\count\z@<0\multiply\count\z@\m@ne\borb_{\the\count\z@}
    \else\oorb{\the\count\z@}\fi\fi
  \futurelet\@next\sep@}%
 \def\o@{\afterassignment\o@@\count\z@=}%
 \def\sep@##1{\ifx\@next]\let\o@)\else
  \ifx\@next,\orbitcup\else\ifx\@next+\orbitplus\fi\fi\fi\o@}%
 \o@}
\makeatother
\def\BB{\Bnd}
\def\BBi{\Bnd'}
\def\BBii{\Bnd''}
\def\BBs{\Bnd^*}
\pdef\B#1[{\BB_{#1}(\orbitlist}
\pdef\Bi#1[{\BBi_{#1}(\orbitlist}
\pdef\Bii#1[{\BBii_{#1}(\orbitlist}
\pdef\Bs#1[{\hyperlink{@-B*}{\BBs_{#1}}(\orbitlist}
\pdef\BL#1[#2;{\BBii_{#1}(#2;\orbitlist}
\let\ddefect\delta
\pdef\defect[#1;{\ddefect(#1;\orbitlist}

\def\atbeginorbits{\iffullversion\else\table[tp]
 \global\advance\tableindex1\relax
 \caption{The lattice \lname{\thelattice}%
  \ifnum\tableindex>1\relax\space\rm(continued)\fi}\label{tab.\thelattice.\the\tableindex}%
 \hrule height0pt\relax\vskip-\abovedisplayskip\small\fi}
\def\atendorbits{\iffullversion\else\endtable\fi}

\def\ssingle{\bar{\frak{s}}_*}
\def\single{\hyperref[s.single]\ssingle}

\let\orbitcup,
\def\orbitplus{\mathbin{\hyperref[rem.comma]{\scriptstyle\cup}}}

\def\mega{\frak m}
\def\ssbnd{\frak b}
\def\sbnd{\hyperlink{@-bnd}{\ssbnd}}
\def\ggeom{\frak B}
\def\geom{\hyperlink{@-geom}{\ggeom}}
\def\patL#1{\hyperlink{@-patL}{\pat_{#1}}}
\def\pperm{\bar\tau}
\def\perm{\hyperlink{@-tau}{\pperm}}
\def\hhyp{\operatorname{hp}}
\def\hyp#1{\hyperlink{@-hyp}{\hhyp_{#1}}}
\def\HHyp{\operatorname{K}}
\def\Hyp{\hyperlink{@-Hyp}{\HHyp}}
\newcommand\OGL[1][\fL]{\hyperlink{@-OG}{\OG}(#1)}
\newcommand\RGL[1][\fL]{\hyperlink{@-OG}{\RG}(#1)}

\def\gchecked{\rlap{\quad\color{green}\checkmark\!\!\!\checkmark}}
\def\gchecked{}

\def\ltitle{\subsection{The lattice \lname{\thelattice}}}
\def\ltext{%
 There are \lcount{\thelattice} isomorphism classes of square $12$
 vectors $\hbar\in\text{\lname{\thelattice}}$,
 with the maximal na\"{\i}ve bound $\bnd(\Orb)=\lmax{\thelattice}$%
 \iffullversion.\else\ (see \alltables).\fi}
\def\latticetext{\iffullversion\ltitle\ltext\else
 \ifnum\lmax{\thelattice}>\mincount\ltitle\ltext\fi\fi
}

\def\hitemii#1{\item[\hlink{\thelattice}{#1}:]\label{s.\thelattice-#1}}
\let\hitem\hitemi

\def\specvec#1{\bar{#1}}

\def\vv#1{\specvec1_{#1}}

\def\cw#1{[#1]}

\def\configform#1{\##1}
\newcommand\config[2][\thelattice]{\hyperlink{#1-#2}{$\configform{#2}$}}

\let\class=c
\def\Classesi{\Cal{C}'}
\def\Classesii{\Cal{C}''}

\def\oind#1{\obj{#1}{#1}}
\def\(#1,#2){(\oind{#1},\oind{#2})}

\def\CH{\Cal H}
\def\CL{\Cal L}
\def\bL{\bold L}
\def\ppos{\Cal H}
\def\pos{\hyperlink{@-pos}{\ppos}}
\def\Golay#1{\CC_{#1}}
\def\units{^\times}
\def\sat{\operatorname{sat}}
\def\val{\operatorname{val}}

\let\graph\fL
\def\ts{\tilde{\frak{s}}}
\def\spnc{\spn^\frak{s}}

\def\loccit{\latin{loc\PERIOD~cit}}

\def\sdif{\mathbin\vartriangle}
\def\CA{\Cal A}

\def\minitab#1{\vcenter\bgroup\rm
 \def\-##1{\setbox0\hbox{$00$}\hbox to\wd0{\hss$##1$\hss}}%
 \let\\\cr
 \halign\bgroup\strut\hss$##$&&#1\hss$##$\hss\cr}
\def\endminitab{\crcr\egroup\egroup}

\allowdisplaybreaks

\title{Conics in sextic {\4$K3$-}surfaces in $\mathbb{P}^4$}

\author{Alex Degtyarev}

\address{%
Department of Mathematics\\
Bilkent University\\
06800 Ankara, TURKEY}

\email{
degt@fen.bilkent.edu.tr}

\thanks{%
The author was partially supported by the T\"{U}B\DOTaccent{I}TAK grant 118F413%
}



\begin{document}

\begin{abstract}
We prove that the maximal number of conics in a smooth sextic {\4$K3$-}surface
$X\subset\mathbb{P}^4$ is 285, whereas the maximal number of real conics in a real
sextic is 261. In both extremal configurations, all conics are irreducible.
\end{abstract}

\maketitle


\section{Introduction}\label{Intro}

All algebraic varieties considered in the paper are over~$\C$.
{\4A \emph{sextic \rom(surface\rom)} is a degree six $K3$-surface
in~$\Cp4$; it is the intersection of a cubic and a quadric.}

\subsection{The principal result}
Counting rational curves in projective surfaces is a long standing classical
problem in algebraic geometry. Each smooth cubic surface $X\subset\Cp3$
contains exactly $27$ lines (A.~Cayley~\cite{Cayley:cubics}, 1849), whereas a
typical quartic surface has no lines at all, and the question changes to
\emph{bounding} the number of lines that a special surface may have.
In 1882, F.~Schur~\cite{Schur:quartics} constructed a quartic with $64$
lines, but the problem
was not settled until 1943, when B.~Segre~\cite{Segre}
proved that $64$ is indeed the maximum. A minor gap in Segre's proof was
discovered and
bridged by S.~Rams and M.~Sch\"{u}tt~\cite{rams.schuett}, and a partial
classification of large configurations of lines was given by A.~Degtyarev,
I.~Itenberg, and A.~S.~Sert\"{o}z~\cite{DIS}.

Little is known about surfaces $X\subset\Cp3$ of higher degree: in
spite of considerable efforts (see recent papers~\cite{Elkies}
or~\cite{rams.schuett:quintics} and references therein), the known
bounds are still too far from the known examples. On the other hand,
smooth quartics being $K3$-surfaces, one can change the paradigm and consider
other polarizations, \ie, smooth  $K3$-surfaces
$X\subset\Cp{d+1}$ of degree~$2d$. In this direction, we do obtain sharp
upper bounds and a plethora of examples, see~\cite{degt:lines,degt:sextics}.

Even less is known about rational curves of other degree, even in
$K3$-surfaces $X\subset\Cp{d+1}$.
For conics in quartic surfaces, the best examples have
$352$
or $432$ conics (see~\cite{Barth.Bauer:conics,Bauer:conics}),
whereas the best known upper
bound is $5016$, see~\cite{Bauer:conics} with a reference to S.~A.~Str{\o}mme.
(There
also are some asymptotic bounds as $d\to\infty$,
see \cite{Miyaoka,rams.schutt:24curves}.)
In this paper, we suggest to attack the conic problem using
a modification of an approach that already proved useful in two
similar counting problems (\cf.~\cite{degt:sextics,DIO}). We expect it to work
for the ``interesting'' degrees $2\le2d\le8$, where the surfaces are complete
intersections, and we test it in the case $2d=6$,
sextic surfaces in~$\Cp4$,
mainly
because we can partially reuse some results of~\cite{DIO} in the computation.
It is worth emphasizing that, in our approach,
we intentionally cut off the lines and work with an abstract combinatorial
configuration (dual adjacency graph) of conics,
counting both irreducible and reducible ones.
Lines are recovered at the
very end, when the lattice is embedded to $H_2(X)$,
and it is at this point that reducible
conics can be told apart.
The combinatorial configuration does
\emph{not} always know which of the conics are reducible, see the remark
after \autoref{th.main} and Examples~\ref{ex.249}, \ref{ex.231}.

The
principal result of the paper is the following sharp upper bound
on the number of conics in a smooth sextic surface $X\subset\Cp4$,
including a partial classification of large combinatorial configurations
of conics.
Surprisingly, the
extremal surfaces contain no lines, so that all their conics are irreducible.

\theorem[see \autoref{proof.main}]\label{th.main}
Up to the projective group, there are but three smooth sextic
{\4$K3$-}surfaces
$X\subset\Cp4$ with more than $260$ conics \rom(irreducible or
reducible\rom)\rom:
\roster
\item\label{main.285}
one real surface $\quartic1$ with $285$ conics, all irreducible,
see \eqref{eq.quartic.1} below\rom;
\item\label{main.261}
two real surfaces $\quartic2$, $\quartic3$ sharing the same
combinatorial
configuration of $261$ conics, all irreducible,
see \eqref{eq.quartic.2}, \eqref{eq.quartic.3} below, respectively.
\endroster
The surface~$\quartic3$ admits a real form in which all
conics are real \rom(moreover, have real points\rom),
thus maximizing the number of \emph{real} conics in a
\emph{real} sextic surface.
\endtheorem

{\4In view of \cite{Schutt.singular.K3,Shimada:conjugate},
the polarized surfaces
$\quartic2$ and $\quartic3$ must be Galois conjugate;
that would explain the fact that they share the same configuration of
conics.
}

Conjecturally, the next value taken by the number of conics in a smooth
sextic surface $X\subset\Cp4$ is $249$, see \autoref{ex.249}. The known
combinatorial configuration of $249$ conics is realized by four distinct surfaces;
in two of these surfaces, all conics are irreducible, whereas in the two others, there
are $60$ irreducible conics and $189$ reducible ones, composed of the $42$
lines contained in the surface.
In fact, $189$ is the maximal possible number of reducible conics, see
\autoref{rem.189}.

\subsection{Contents of the paper}
In \autoref{S.reduction}, we introduce the arithmetical reduction of the
problem, consisting in modifying the N\'{e}ron--Severi lattice of the surface
and replanting the result to an appropriate Niemeier lattice. With future
applications in mind, we treat
arbitrary smooth $2d$-polarized $K3$-surfaces $X\to\Cp{d+1}$. In
\autoref{S.computation}, we describe the principal technical tools used in the
computation. We switch to degree $2d=6$, even though most tools should
apply to other degrees as well. \S\S\,\ref{S.few}--\ref{S.Leech} deal with
the computation \latin{per se}, which was done using \GAP~\cite{GAP4}; we
outline just enough details for the interested reader to be able to reproduce
the result. Finally, in \autoref{S.proofs} we combine our findings to present
a formal proof of \autoref{th.main}. We also discuss a few interesting
examples of large configurations of conics which were discovered in our study
of individual Niemeier lattices.

\subsection{Acknowledgements}
I am grateful to S{\l}awomir Rams for introducing me to the subject and
bringing to my attention
the few known results.
{\4I would also like to thank J\'{a}nos Koll\'{a}r for a few helpful remarks and
the anonymous referee of this paper for a number of suggestions improving the
text.}


\section{The reduction}\label{S.reduction}

Throughout the paper, we fix the notation~$\bL$ for a
fixed even unimodular lattice of
signature $(3,19)$; in other words, $\bL\cong-2\bE_8\oplus3\bU$ is the
$2$-homology of a $K3$-surface.
All lattices considered are \emph{even}, and a \emph{root} is a vector
of square $\pm2$.

\subsection{Conics as homology classes}\label{s.homology}
Let $\Gf\:X\into\Cp{d+1}$ be a smooth $2d$-polarized $K3$-surface.
(If $d=1$, the ``smoothness'' refers to
the ramification locus of~$\Gf$, which is two-to-one; otherwise, we assume
$\Gf$ injective.) The map~$\Gf$ gives rise to the triple
\[*
\bL\cong H_2(X)\supset\NS(X)\ni h_X:=\Gf^*\CO_{\Cp{d+1}}(1),\quad h_X^2=2d,
\]
called the \emph{homological type} of $X\into\Cp{d+1}$. The \emph{period}
$\Go_X\in H_2(X;\C)$ (the class of a holomorphic
$2$-form on~$X$) defines an orientation on
$\R h\oplus\R(\Re\Go_X)\oplus\R(\Im\Go_X)$ and, hence, on any other
positive definite $3$-space in $H_2(X;\R)$; this coherent choice of
orientations is called the \emph{canonical orientation} of the homological
type.
Letting $(\bL\supset\NS\ni h):=(H_2(X)\supset\NS(X)\ni h_X)$,
the homological type has the following properties:
\roster
\item\label{i.hyperbolic}
the lattice $\NS$ is \emph{hyperbolic}, \ie, $\Gs_+\NS=1$;
\item\label{i.primitive}
the sublattice $\NS\subset\bL$ is \emph{primitive}, \ie,
$\Tors(\bL/\!\NS)=0$;
\item\label{i.root}
there is no vector $e\in\NS$ such that $e^2=-2$ and $e\cdot h=0$\rom;
\item\label{i.isotropic}
there is no vector $e\in\NS$ such that $e^2=0$ and $e\cdot h=2$\rom;
\item\label{i.Veronese}
{\4if $d=4$, then $h\notin2\NS$.}
\endroster
{\4
Here, \iref{i.isotropic} implies that there is no $e\in\NS$ such that $e^2=0$
and $e\cdot h=1$ and, hence, the linear system~$\ls|h|$ is fixed point
free, see \cite{Nikulin:Weil,Saint-Donat}. Then, by~\cite{Saint-Donat},
conditions~\iref{i.isotropic} and~\iref{i.Veronese} are equivalent to the
statement that the map $X\to\Cp{d+1}$
given by $\ls|h|$ is birational onto its image. (Strictly speaking,
\cite{Saint-Donat} uses a genus~$1$ curve rather than class $e\in\NS$; however,
assuming~\iref{i.root}, one of $\pm e$ is necessarily nef.)}

Due to the global Torelli theorem for
$K3$-surfaces~\cite{Pjatecki-Shapiro.Shafarevich}, surjectivity of the period
map~\cite{Kulikov:periods}, and Saint-Donat's results on projective models of
$K3$-surfaces~\cite{Saint-Donat},
the converse also holds.
{\4The following statement is well known, and various versions thereof appear
in virtually any paper on the subject.}

\theorem[{\cf. \cite[Theorem 3.11]{DIS} or \cite[Theorem 7.3]{degt:singular.K3}}]\label{th.existence}
An oriented abstract triple $\bL\supset\NS\ni h$, $h^2=2d$, is
isomorphic to the
oriented
homological type of a smooth $2d$-polarized
$K3$-surface~$X$  if and only if $\bL\supset\NS\ni h$
has properties~\iref{i.hyperbolic}--\iref{i.Veronese} above.
In this case, $X$ is a generic member of a connected
$(20-\rank\NS)$-parameter \rom(modulo the
group
$\PGL(d+2,\C)$\rom)
family of $2d$-polarized $K3$-surfaces.
\done
\endtheorem

Given a line $L\subset X$ or a conic $C\subset X$ (not necessarily
irreducible), the classes $l:=[L]\in\NS(X)$ and $c:=[C]\in\NS(X)$ have the
property that
\[
l^2=c^2=-2,\qquad l\cdot h=1,\quad c\cdot h=2.
\label{eq.lch}
\]
In particular, since the self-intersection is negative, a line or a conic is
unique in its class and can be identified with the latter.
By the Riemann--Roch theorem, conversely, any class $l$ or $c\in\NS(X)$
satisfying~\eqref{eq.lch} is represented by a unique line or conic,
respectively.
It follows from the description of the nef cone of~$X$
(see, \eg, \cite[Corollary~8.2.11]{Huybrechts}) that \emph{irreducible}
conics are characterized by the property
\[
\text{a conic $c\in\NS(X)$ is irreducible iff $c\cdot l\ge0$ for each line
 $l\in\NS(X)$}.
\label{eq.irreducible}
\]

\subsection{The lattice $S(\NS,h)$}\label{s.S-lattice}
A \emph{$2d$-polarized lattice} is a lattice $\NS$ equipped with a
distinguished vector $h\in\NS$ of square $2d>0$.
We say that $\NS\ni h$ is of
\emph{type~$\I$} if $h\in2\NS\dual$; otherwise, $\NS\ni h$ is said to
be of \emph{type~$\II$}. The \emph{graph of lines} ($n=1$), \emph{conics}
($n=2$) \etc. of a
\emph{hyperbolic} polarized lattice $\NS\ni h$ is the graph
\[*
\Fn_n(\NS,h):=\bigl\{c\in\NS\bigm|
 \text{$c^2=-2$ and $c\cdot h=n$}\bigr\},
\]
where two vertices $c_1$, $c_2$ are connected by an edge of multiplicity
$c_1\cdot c_2$.
Since $\NS$ is assumed hyperbolic, each set $\Fn_n$ is finite;
the vertices of $\Fn_1$ and $\Fn_2$ are called \emph{lines}
and \emph{conics}, respectively.
{\4Obviously, if there is at least one conic, the extension
$\NS\supset h^\perp\oplus\Z h$ is of index $d$ (type~$\I$) or $2d$ (type~$\II$),
and it is these two cases that we are considering below.
If there is at least one line, the extension is of index~$2d$
and, hence, $\NS\ni h$ is necessarily of type~$\II$.}

Fix a hyperbolic polarized lattice $\NS\ni h$ and assume that
$\Fn_2(\NS,h)\ne\varnothing$.
Consider the orthogonal complement $h^\perp$ and
project each conic $c\in\NS$ to
\[*
c\mapsto c':=c-d\1h\in(h^\perp)\dual,\quad c^{\prime2}=-2(d+1)d\1.
\]
Then, consider the lattice
\[*
S:=-(h^\perp\oplus\Z\hbar)_d^\sim,\quad \hbar^2=2d(1-d);
\]
here, $^\sim$ stands for the cyclic index~$d$ extension generated by
$c'+d\1\hbar$, where $c'$ is the image of any conic $c\in\NS$. (In the
exceptional case $d=1$, we merely let $\hbar=0$ and $S:=-h^\perp$.)
The result is a polarized positive definite lattice
\[*
S:=S(\NS,h)\ni\hbar,\quad \hbar^2=2d(d-1),\quad \hbar\in2(d-1)S\dual.
\]
Furthermore, there are bijections between the following sets:
\roster
\item\label{bijection.conics}
conics $c\in\NS$
$\longleftrightarrow$
vectors $l\in S$ such that $l^2=4$,
$l\cdot\hbar=2(d-1)$;
\item\label{bijection.exceptional}
\emph{exceptional divisors} $e\in h^\perp$, $e^2=-2$
$\longleftrightarrow$
roots $r\in\hbar^\perp\subset S$,
\cf. \autoref{s.homology}\iref{i.root};
\item\label{bijection.isotropic}
\emph{$2$-isotropic vectors} $e\in\NS$, $e^2=0$, $e\cdot h=2$
$\longleftrightarrow$
roots $r\in S$ such that $r\cdot\hbar=2(d-1)$,
\cf. \autoref{s.homology}\iref{i.isotropic}.
\endroster
{\4(The bijection in~\iref{bijection.exceptional} is the ``identity'', and
in~\iref{bijection.conics},~\iref{bijection.isotropic} it is
$v\mapsto v-d\1h+d\1\hbar$.)}
In particular,
we have the following statement.

\lemma\label{lem.root.free}
A polarized hyperbolic lattice $\NS\ni h$ has no vectors as in
\autoref{s.homology}\iref{i.root} or \iref{i.isotropic} if and only if
the corresponding lattice $S(\NS,h)$ is root free.
\endlemma

\proof
The statement is obvious in the case $d=1$, where all three sets (exceptional
divisors, $2$-isotropic vectors, and roots in~$S$) are in a bijection with
one another.

If $d\ge2$, since $S$ is positive definite, we have
$r\cdot\hbar\le2\sqrt{d(d-1)}<4(d-1)$ for each root $r\in S$. Since also
$\hbar\in2(d-1)S\dual$, the product $r\cdot\hbar$ can take but three values
$0,\pm2(d-1)$, and the statement follows from \iref{bijection.exceptional}
and \iref{bijection.isotropic} above.
\endproof

\subsection{Intersection of conics}\label{s.intr}
Let $S\ni\hbar$ be a positive definite $2d(d-1)$-polarized lattice.
Assuming $S$ root free, we define its \emph{Fano graph}
\[*
\Fn(S,\hbar):=\bigl\{l\in S\bigm|
 \text{$l^2=4$ and $l\cdot\hbar=2(d-1)$}\bigr\},
\]
where two vertices $l_1$, $l_2$ are connected by an edge of multiplicity
$2-(l_1\cdot l_2)$. This graph appears in a number of geometric problems
(\cf.~\cite{degt:sextics,DIO}), yielding a number of names for its vertices; in
this paper, in view of the canonical graph isomorphism
$\Fn(S,\hbar)=\Fn_2(\NS,h)$ constructed in \autoref{s.S-lattice}, we are
compelled to call them \emph{conics}.
However, we retain the notation~$l$ for a ``typical'' vertex.

In view of the next lemma (reflecting the obvious geometry of a pair of
conics in a $K3$-surface),
all edges of $\Fn(S,\hbar)$ have non-negative multiplicity.

\lemma\label{lem.intr}
If $S$ as above is root free and $\hbar\in S$ is a primitive vector, then,
for any pair $l_1\ne l_2\in\Fn(S,\hbar)$, one has
\roster*
\item
if $d=1$, then $l_1\cdot l_2\in\{-4,-2,-1,0,1,2\}$\rom;
\item
if $d=2$, then $l_1\cdot l_2\in\{-2,0,1,2\}$\rom;
\item
in all other cases, $l_1\cdot l_2\in\{0,1,2\}$.
\endroster
One has $l_1\cdot l_2=-4$ \rom(if $d=1$\rom) or $-2$ \rom(if $d=2$\rom) if
and only if $l_1+l_2=\hbar$.
\endlemma

\proof
Since $S$ is positive definite, one has $\ls|l_1\cdot l_2|\le4$. The
possibility $l_1\cdot l_2=\pm3$ is ruled out since otherwise $l_1\mp l_2$
would be a root. If $d=1$, there are no further restrictions. Otherwise, the
obvious inequality $\det(\Z\hbar+\Z l_1+\Z l_2)\ge0$ implies that
{\4$l_1\cdot l_2\ge-4d\1\ge-2$}, and there remains to rule out the possibility
$l_1\cdot l_2=-1$ for the values $d=2,3,4$.

If $d=2$ or~$3$, then $\hbar-(l_1+l_2)$ is a root.
If $d=4$, then {\4$\hbar-2(l_1+l_2)$ annihilates all three vectors and,
hence,}
$\hbar=2(l_1+l_2)$ is not primitive.
{\4(Recall that any sublattice of a definite lattice is nondegenerate.)}
\endproof

\subsection{The discriminant \pdfstr{discr S}{$\discr S$}}\label{s.discr.S}
Recall that the \emph{discriminant} of a nondegenerate even lattice~$L$ is
the finite abelian group
$\CL:=\discr L:=L\dual\!/L$
equipped with the nondegenerate quadratic form
\[*
\CL\to\Q/2\Z,\quad x\bmod L\mapsto x^2\bmod2\Z.
\]
(The associated symmetric bilinear form takes values in $\Q/\Z$. For
details, see~\cite{Nikulin:forms}, where the notation is $q_L$.)
For a prime~$p$,
we abbreviate
\[*
\CL_p:=\discr_pL:=\CL\otimes\Z_p\quad\text{and}\quad
\ell_p(\CL):=\ell(\CL_p),
\]
where $\ell(A)$ stands for the minimal number of
generators of an abelian group~$A$. The determinant of
$\CL_p$ (see~\cite{Nikulin:forms}) is denoted by $\det_p\CL=\det_p\CL_p$;
it has the form $\Ge_p/\ls|\CL_p|$, where the class
$\Ge_p\in\Z_p\units/(\Z_p\units)^2$ is well defined (and used) unless $p=2$ and
$\CL_2$ is \emph{odd}, \ie, contains an order~$2$ vector~$x$ with
$x^2\ne0\bmod\Z$.
A vector $v\in\CL_2$ is called \emph{characteristic} if
$x^2=x\cdot v\bmod\Z$ for any \emph{order~$2$} vector $x\in\CL_2$.

In this section, we describe the relation between
\[*
\CS_p:=\discr_pS,\quad
\CN_p:=\discr_p{\NS},\quad\text{and}\quad
\CH_p:=\discr_ph^\perp
\]
for each prime~$p$.
As above, we assume that {\em there is at least one conic $c\in\NS$}.
To begin with, we observe that
\[
\ls|\CS|=d'\ls|\CN|\bmod(\Q\units)^2,
\label{eq.|S|}
\]
where $d'=2$ if $d=1$ and $d'=d-1$ otherwise.

Given~$p$,
we denote by $q:=p^r$ the maximal power of~$p$ that divides $2d(d-1)$; in the
special case $d=1$, we let $q=1$ for all $p\ne2$ and $q=2$ for $p=2$.


\subsubsection{The case of $p$ prime to $2(d-1)$ \rom(including $d=1$ and $p\ne2$\rom)}\label{ss.p|d}
If $q=1$, we have an obvious canonical isomorphism
\[
-\CS_p=\CN_p.
\label{eq.S=N}
\]
The same holds for $q>1$.
Indeed, put $d=q\bar{q}$ and let $c'\in(h^\perp)\dual$ be the projection
of a conic $c\in\NS$.
Then, it is obvious from the construction and Nikulin's
theory~\cite{Nikulin:forms} applied to the finite index extensions
${\NS}\supset h^\perp\oplus\Z h$ and
${-S}\supset h^\perp\oplus\Z\hbar$
that
the class
$(\bar{q}c'\bmod h^\perp)\in\CH_p$, which has order~$q$ and square
$-2\bar{q}q\1\bmod2\Z$,
generates a nondegenerate subgroup (hence, orthogonal summand)
$\Z/q\subset\CH_p$
and there are canonical isomorphisms
$\CN_p=\bigl[(\bar{q}c')^\perp\subset\CH_p\bigr]=-\CS_p$.

\subsubsection{The case of $q>1$ and $p$ prime to $2d$}\label{ss.p|(d-1)}
Let $d-1=\bar{q}q$. Then, $d=\bar{q}q+1$ and
there is an obvious canonical isomorphism
\[*
-\CS_p=\CN_p\oplus\bigl\<q\1\hbar\bigr\>,\qquad
\bigl(q\1\hbar\bigr)^2=-2d\bar{q}q\1\bmod2\Z.
\]
In particular, in view of~\eqref{eq.|S|},
\[
\ell_p(\CS)=\ell_p(\CN)+1,\qquad
\ls|\CS|\det\nolimits_p(-\CS)=-2\ls|\CN|\det\nolimits_p\CN\bmod(\Z_p\units)^2.
\label{eq.det}
\]

\subsubsection{The case of $p=2$ and $d>1$ odd}\label{ss.d.odd}
Letting $2(d-1)=q\bar{q}$, we have
\[*
-\CS_2=
\CH_2\oplus\bigl<q\1\hbar\bigr>,
\qquad
(q\1\hbar)^2=-d\bar{q}q\1\bmod2\Z.
\]
If $\NS$ is of type~$\I$, the class \smash{$\frac12h\bmod{\NS}$}
of order~$2$ and
square $\frac12d\bmod2\Z$ generates an orthogonal summand $\Z/2$ in $\CN_2$ and
we have
\[*
\CH_2=\bigl\<\tfrac12h\bigr\>^\perp\subset\CN_2,\quad\text{or}\quad
\CN_2=\CH_2\oplus\bigl\<\tfrac12h\bigr\>.
\]
It follows that
\[*
\ell_2(\CS)=\ell_2(\CN)\quad\text{and}\quad
\text{$\CN_2$ is odd}.
\]
Otherwise, if $\NS$ is of type~$\II$, then
there is a class $\kappa\in\CH_2$ of order~$2$ and square $-\frac12d\bmod2\Z$
such that
\[*
\CN_2=\<\kappa\>^\perp\subset\CH_2.
\]
Hence,
\[*
\ell_2(\CS)=\ell_2(\CN)+2\quad\text{and}\quad
\text{$\CS_2$ is odd}.
\]

\subsubsection{The case of $p=2$ and $d=1$}\label{ss.d=1}
This case is straightforward: $\CS_2=\CH_2$ and
\[*
\CN_2=\CH_2\oplus\bigl<\tfrac12h\bigr>\quad\text{(type $\I$)}
\quad\text{or}\quad
\CH_2=\CN_2\oplus\<\kappa\>\quad\text{(type $\II$)},
\]
where $\kappa$ has order~$2$ and $\kappa^2=-\frac12\bmod2\Z$.

\subsubsection{The case of $p=2$ and $d$ even}\label{ss.d.even}
Denote $2d=q\bar{q}$ and let $\tilde\kappa\in\CH_2$ be the class of
the vector $\bar{q}c-2q\1h$, where $c\in\NS$ is any conic. Then
\[*
\CS_2=\bigl<\tilde\kappa+2q\1\hbar\bigr>^\perp\!/\bigl<\tilde\kappa+2q\1\hbar\bigr>
\quad\text{in}\quad\CN_2\oplus\bigl<q\1\hbar\bigr>.
\]
If $\NS$ is of type~$\I$, then also
\[*
\CN_2=\bigl<\tilde\kappa+2q\1h\bigr>^\perp\!/\bigl<\tilde\kappa+2q\1h\bigr>
\quad\text{in}\quad\CN_2\oplus\bigl<q\1h\bigr>
\]
and, trying one-by-one the few possibilities for $\tilde\kappa\in\CH_2$
(\cf.~\cite{Cisem:real}), we conclude that
\[*
\CN_2\cong\CS_2\ \text{as groups; hence},\ \ell(\CN_2)=\ell(\CS_2).
\]
(There is no universal relation between the determinants of
$\CN_2$ and $\CS_2$ or even their parity.) If $\NS$ is of type~$\II$, then
there is a class $\kappa\in\CH_2$ such that $2\kappa=\tilde\kappa$,
$\kappa^2=-(\bar{q}/q)\bmod2\Z$, and
\[*
\CN_2=\<\kappa\>^\perp\subset\CH_2.
\]
Hence,
\[*
-\CS_2=\CN_2\oplus\bigl(\<\Ga\>+\<\Gb\>\bigr),\quad\text{where}\quad
\Ga:=\tfrac12q\kappa,\quad\Gb:=\kappa+q\1\hbar
\]
generate a subgroup $\Z/2\oplus\Z/2\subset\CS_2$ and
\[*
\Ga^2=-\tfrac12d\bmod2\Z,\quad
\Ga\cdot\Gb=\tfrac12\bmod\Z,\quad
\Gb^2=-\tfrac{1}{2}\bmod2\Z.
\]
It follows that
\[*
\text{$\CS_2$ is odd}\quad\text{and}\quad
\ell_2(\CS)=\ell_2(\CN)+2.
\]

\subsection{Embeddings to Niemeier lattices}\label{s.embeddings}
An embedding $S\into N$ of free abelian groups is called
\emph{$p$-primitive} for a prime~$p$ if
$\Tors_p(N/S)=0$. The embedding
is called \emph{primitive} if it is $p$-primitive for
each prime~$p$.

\lemma\label{lem.extension}
Any non-degenerate lattice~$S$ admits a finite index extension $S'\supset S$
such that, for each prime~$p$, either
\roster*
\item
$\ell(\discr_pS')\le3$ and $\discr_pS'$ odd \rom(for $p=2$ only\rom), or
\item
$\ell(\discr_pS')\le2$.
\endroster
\endlemma

\proof
It follows immediately from the classification of finite quadratic forms
(see, \eg, \cite{Nikulin:forms}) that
any form on a $p$-group~$\CS$ that does not satisfy the bounds in the statement has
a non-trivial isotropic subgroup $\CK\subset\CS$; then, for the extension
$S'\supset S$ defined by~$\CK$, one has
\[*
\ls|\discr_pS'|=\ls|\CK^\perp\!/\CK|<\ls|\CS|.
\]
Arguing by induction, one reduces the group $\discr_pS'$ to the
bounds stated.
\endproof

\proposition\label{prop.embedding}
Assume that $\NS\ni h$ is a primitive sublattice of~$\bL$ and denote
$S:=S(\NS,h)$. Then\rom:
\roster
\item\label{embedding.S}
$S$ admits a primitive isometric embedding to a Niemeier lattice~$\bN$\rom;
\item\label{embedding.S+A}
the lattice $S\oplus\bA_1$ admits an embedding to a Niemeier
lattice $\bN$ which is $p$-primitive for all \emph{odd} primes~$p$ other than
$p=(4k-1)\divides|(d-1)$.
\endroster
If $\NS\ni h$ is of type~$\I$, then
\roster[\lastitem]
\item\label{embedding.S+A.2.I}
there is an embedding $S\oplus\bA_1\into\bN$ as in \iref{embedding.S+A} which
is also $2$-primitive.
\endroster
\endproposition

\proof
The statement follows from \cite[Theorem 1.12.2]{Nikulin:forms},
the computation of $\CS$ in \autoref{s.discr.S},
and the fact that $\rank\bN=24=\rank\bL+2$.
Indeed, for all primes~$p$,
\[*
\ell(\CN_p)+\rank{\NS}\le22,\quad\text{and, hence},\quad
\ell(\CS_p)+\rank{S}\le24,
\]
with the latter inequality strict unless
{\4$p=2$ and}
$\CS_p$ is odd. Hence, there is no
obstruction to a primitive embedding $S\into\bN$.
Furthermore, letting $\bar S:=S\oplus\bA_1$, for $p>2$ we have
\[*
\ell(\bar\CS_p)+\rank{\bar S}\le24,\quad
\ls|\bar\CS|\det\nolimits_p(-\bar\CS)=-\ls|\CN|\det\nolimits_p\CN\bmod(\Z_p\units)^2,
\]
\cf.~\eqref{eq.det}. The equality $\ell(\bar\CS_p)+\rank{\bar S}=24$ may hold
only if $d>1$, $p\divides|(d-1)$, and
\[*
\ell(\CN_p)+\rank{\NS}=22,\quad\text{hence}\quad
\ls|\CN|\det\nolimits_p\CN\in(\Z_p\units)^2.
\]
This creates an obstruction to a primitive embedding $\bar S\into\bN$ only if
$-1\notin(\Z_p\units)^2$, \ie, $p=3\bmod4$. In this case, since
$\ell(\bar\CS_p)=\rank\bN-\rank\bar S\ge3$, we use \autoref{lem.extension} and
pass to a finite index extension $\bar S'\supset\bar S$ with
$\ell(\bar\CS'_p)\le2$. Similarly, if $p=2$ and $\NS$ is of type~$\II$, it
may happen that
\[*
\ell(\bar\CS_2)>\rank\bN-\rank\bar S\ge3,
\]
in which case \autoref{lem.extension} needs to be used for $p=2$.
\endproof

\section{The computation}\label{S.computation}

From now on, we confine ourselves to the case $d=3$,
{\4\ie, we consider conics in smooth sextic $K3$-surfaces $X\subset\Cp4$.}

\subsection{Admissible and geometric sets}\label{s.admissible}
Following the idea outlined in \autoref{s.S-lattice}, we are interested in
finding root-free $12$-polarized lattices $S\ni\hbar$, admitting an embedding to a
Niemeier lattice~$\bN$ (and satisfying a few extra conditions), with a large
Fano graph $\ls|\Fn(S,\hbar)|>\the\minmin$.
We choose to employ \autoref{prop.embedding}\iref{embedding.S}, insisting
that the embedding $S\into\bN$ should be primitive. Since we can also assume
that $S\otimes\Q$ is generated by~$\hbar$ and the conics, we can build this
lattice directly inside~$\bN$.

More precisely, we fix a Niemeier lattice~$\bN$ and a class $\hbar\in\bN$,
$\hbar^2=12$, and consider the set
\[*
\fF:=\fF(\bN,\hbar):=\bigr\{l\in\bN\bigm|
 \text{$l^2=4$ and $l\cdot\hbar=2(d-1)$}\bigr\}
\]
of ``prospective conics.'' (Since $\bN$ is not assumed {\4root} free, we do not
regard this set as a graph and use a different notation.) For a subset
$\fL\subset\fF$, we introduce the following terminology and notation:
\roster*
\item
the \emph{span} $\spn\fL:=(\Q\fL+\Q\hbar)\cap\bN$ and
\emph{rank} $\rank\fL:=\rank\spn\fL$;
\item
occasionally, we use the \emph{integral span}
$\spn_\Z\fL:=(\Z\fL+\Z\hbar)\subset\bN$;
\item
a set $\fL$ is \emph{generated} by a subset $\fL'\subset\fL$ if
$\spn\fL'=\spn\fL$;
\item
the \emph{saturation} of~$\fL$ is the set $\sat\fL:=\fF\cap\spn\fL\subset\fF$;
\item
a set $\fL$ is \emph{saturated} if $\fL=\sat\fL$;
\item
similarly, \via\ $\spn_\Z$, we define \emph{$\Z$-saturation} $\sat_\Z$ and
\emph{$\Z$-saturated sets};
\item
a saturated set~$\fL$ is \emph{admissible} if $\spn\fL$ is root free
and $\hbar\in4(\spn\fL)\dual$.
\endroster
Recall that the last condition implies that
\roster*
\item
$l_1\cdot l_2\in\{0,1,2\}$ for any two conics $l_1\ne l_2\in\fL$,
\endroster
see \autoref{lem.intr}; technically, this condition is usually checked first,
and it is
this condition alone that rules out most pairs $\bN\ni\hbar$.


The lattice $\spn\fL\ni\hbar$ plays the r\^{o}le of
$S\ni\hbar$ in \autoref{s.S-lattice}.
{\4We denote by $\hyp\kappa(\fL)$ the result of applying to
$S:=\spn\fL\ni\hbar$ the inverse construction of \autoref{s.S-lattice}.
It}
depends on an extra parameter $\kappa=0$
(type~$\I$) or $\kappa\in\Hyp(\fL)$:
\[
\anchor[@]{Hyp}\HHyp(\fL):=\bigl\{\kappa\in\discr_2(\spn\fL)\bigm|
 \text{$\kappa^2=\tfrac32\bmod2\Z$,
 $\kappa\cdot\hbar=0\bmod4$,
 $2\kappa=0$}
 \bigr\}
\label{eq.Hyp}
\]
(type~$\II$), \cf.~\autoref{ss.d.odd}.
{\4The
construction of $\hyp0(\fL)$ is straightforward, \cf.
\autoref{s.S-lattice}:
it is the cyclic index~$d$ extension of $\hbar^\perp\oplus\Z h$ containing
some (equivalently, any) vector of the form $l-d\1\hbar+d\1h$, $l\in\fL$.
Then,}
$\anchor[@]{hyp}{\hhyp_\kappa(\fL)}\supset\hhyp_0(\fL)$ is the index~$2$
extension defined by $\kappa+\frac12h$.

An admissible
set~$\fL$ is called
\emph{geometric} if at least one of these lattices $\hyp\kappa(\fL)$,
{\4$\kappa\in\Hyp(\fL)\cup\{0\}$,}
admits a primitive isometric embedding to~$\bL$.
Using Nikulin~\cite{Nikulin:forms}, the
{\4computation in}
\autoref{s.discr.S} can be
recast to the following criterion.

\proposition\label{prop.I}
Let $\fL\subset\fF$ be an admissible set, and denote $\CS:=\discr(\spn\fL)$
and $r:=\rank\fL$. Then,
the type~$\I$ lattice $\hyp0(\fL)$ admits a primitive embedding to~$\bL$
if and only if the following
statements hold\rom:
\roster
\item\label{I.rank}
$r\le20$\rom;
\item\label{I.p}
for each odd prime~$p$, one has $\ell(\CS_p)\le20-r$ and,
if $\ell(\CS_p)=20-r$, the congruence
$\ls|\CS|\det_p(-\CS)=2\bmod(\Z_p\units)^2$ holds,
\cf.~\eqref{eq.|S|} and~\eqref{eq.S=N}\rom;
\item\label{I.2}
$\ell(\CS_2)\le20-r$.
\endroster

The type~$\II$ lattice $\hyp\kappa(\fL)$, $\kappa\in\Hyp(\fL)$,
admits a primitive embedding to~$\bL$ if and only if
statements~\iref{I.rank}, \iref{I.p} above hold and
either \iref{I.2} holds or
\roster[\lastitem]
\item\label{II.2}
$\ell(\CS_2)=22-r$ and
\roster*
\item
$\kappa$ is characteristic in $\CS$ and
$\frac12\ls|\CS|\det_2(\kappa^\perp)=\pm3\bmod(\Z_2\units)^2$ or
\item
$\kappa$ is not characteristic in~$\CS$.
\done
\endroster
\endroster
\endproposition

The following simple observation is crucial: it ensures that, when
constructing a geometric set inductively, adding a few conics at a time, we
can discard immediately all intermediate sets that are not geometric.

\lemma\label{lem.hereditary}
Both ``admissible'' and ``geometric'' are hereditary properties\rom:
if a set $\fL\subset\fF$ is admissible/geometric, so is any saturated subset
$\fL'\subset\fL$.
\done
\endlemma

\subsection{Orbits and bounds}\label{s.orbits}
Fix $\bN\ni\hbar$ and $\fF:=\fF(\bN,\hbar)$ as above and introduce the groups
\roster*
\item
$\OG_\hbar(\bN)$, the stabilizer of~$\hbar$ in $\OG(\bN)$;
\item
$\RG_\hbar(\bN)\subset\OG_\hbar(\bN)$, the stabilizer of~$\hbar$ in
the normal subgroup $\RG(\bN)\subset\OG(\bN)$ generated by reflections
$\tr_r\:x\mapsto x-(x\cdot r)r$ defined by roots
$r\in\bN$;
\item
the \emph{stabilizers}
\anchor[@]{OG}$\OG(\fL)\subset\OG_\hbar(\bN)$ and
\anchor[@]{RG}$\RG(\fL)\subset\RG_\hbar(\bN)$
of a geometric set~$\fL$;
\item
the stabilizer
$\stab\hbar:=\OG_\hbar(\bN)/\!\RG_\hbar(\bN)$.
\endroster
The set~$\fF$ splits into $\OG_\hbar(\bN)$-\emph{orbits},
$\fF=\bigcup_n\borb_n$, and each orbit $\borb_n$ splits into
$\RG_\hbar(\bN)$-orbits, referred to as the \emph{combinatorial orbits}.
The group $\stab\hbar$ acts on
\roster*
\item
the set $\Orb:=\Orb(\bN,\hbar)$ of all combinatorial orbits,
\endroster
and $\borb_n$ can as well be regarded as the orbits of this action.
(To simplify the notation, we abuse the language and sometimes treat
subsets of~$\Orb$ as sets of {\4conics (subsets of~$\fF$)};
conversely, the orbits $\borb_n$ are
often treated as subsets of~$\Orb$ rather than~$\fF$.)

A typical combinatorial orbit~$\orb$ is relatively small, and it is easy to
compute
\roster*
\item
the \emph{geometric intersections}
$\anchor[@]{geom}\ggeom(\orb):=
\bigl\{\fL\cap\orb\bigm|\text{$\fL\subset\fF$ is geometric}\bigr\}$
and
\item
the set
$\anchor[@]{bnd}{\sbnd(\orb)}:=\bigl\{\ls|\fL|\bigm|\fL\in\geom(\orb)\bigr\}$
and the \emph{bound}
$\bnd(\orb):=\max\sbnd(\orb)$.
\endroster
(Certainly, it suffices to consider geometric sets~$\fL$ generated by subsets
of~$\orb$.) Then, given a subset $\Cluster\subset\Orb$, we \emph{define} the
na\"{\i}ve \latin{a priori} bound \via
\roster*
\item
$\bnd(\Cluster):=\sum\bnd(\orb)$, $\orb\in\Cluster$.
Obviously, $\ls|\fL\cap\Cluster|\le\bnd(\Cluster)$ for any geometric set~$\fL$.
\endroster
At this point, we discard the pair $\bN\ni\hbar$ if
$\bnd(\Orb)\le\the\minmin$. (The vast majority of pairs is ruled out
by the rough combinatorial estimates found in~\cite{DIO}.) We define
\roster*
\item
the \emph{defect}
$\ddefect(\fL;\Cluster):=\bnd(\Cluster)-\ls|\fL\cap\Cluster|$ of a geometric
set~$\fL$ and
\item
the set
$\BB_d(\Cluster):=\{\fL\subset\fF\,|\,
 \text{$\fL$ is geometric and $\ddefect(\fL;\Cluster)\le d$}\}/\!\OG_\hbar(\bN)$,
$d\in\NN$.
\endroster

We are interested in large geometric sets satisfying the inequality
\[
\ls|\fL|\ge M:=261.
\label{eq.threshold}
\]
Hence, our ultimate goal is the computation of the set
$\BB_{\bnd(\Orb)-M}(\Orb)$.
Below, without further reference, we use the following simple observation.

\lemma\label{lem.bounds}
Let $d_i\in\NN$ and $\Cluster_i\subset\Orb$, $1\le i\le n$,
be integers and pairwise disjoint subsets such that
\[*
d_1+\ldots+d_n+n>\bnd(\Orb)-M.
\]
\roster
\item
{\4If $\BB_{d_i}(\Cluster_i)=\varnothing$ for all $1\le i\le n$,
then $\BB_{\bnd(\Orb)-M}(\Orb)=\varnothing$.}
\item
{\4If $\BB_{d_i}(\Cluster_i)=\varnothing$ for all $2\le i\le n$ and
the inequality $\ddefect(\fL;\Cluster_2)>d_1+d_2+1$
holds for each set $\fL\in\BB_{d_1}(\Cluster_1)$,
then $\BB_{\bnd(\Orb)-M}(\Orb)=\varnothing$.}
\done
\endroster
\endlemma

In the rest of this section, we outline,
{\4for future reference,}
a few technical aspects concerning
the computation of the
{\4bounds $\bnd(\orb)$ and}
sets $\BB_d(\Cluster)$. Most details are found
in~\cite{degt:sextics,DIO},
{\4and the computation itself is contained in
\autoref{S.few}--\autoref{S.Leech}.}

\subsection{Combinatorial estimates\pdfstr{}{ \rm(see~\cite{DIO})}}\label{s.estimates}
{\4We always start the computation from simple \emph{combinatorial estimates}
$\bnd'(\orb)\ge\bnd(\orb)$ for each combinatorial orbit~$\orb$.
They are computed as explained in~\cite{DIO}. Alternatively, in more details,
this block-by-block techniques is presented in \cite[\S\,4]{degt:sextics}
(see $\bnd(\orb)$ in Eqn.\spacefactor1000~4.3), to
which a few minor adjustments are to be made.
First, one should disregard all statements in~\cite{degt:sextics} involving
the \emph{duality} $l\leftrightarrow l^*$. Second,
equation~$(4.1)$ in~\cite{degt:sextics} is to be replaced with
\[
\text{$l^{\prime2}-l'\cdot l''=0$ (iff $l'=l''$), $2$, $3$, or $4$
 (\cf. \autoref{lem.intr})}
\label{eq.intr}
\]
for any two conics $l',l''$ (or their projections to a block). Finally, in
view of the latter change, the combinatorial Lemmas~4.7 and~4.8
in~~\cite{degt:sextics} are replaced, respectively, with the following
simple rough bounds, which suffice for our purposes.}

\lemma[see~\cite{DIO}]\label{lem.max.set.A}
{\4Consider a finite set~$S$, $\ls|S|=n$, and let $\fS$ be a collection of
subsets $s\subset S$ with the following properties\rom:
\roster
\item\label{sets.m}
all subsets $s\in\fS$ have the same fixed cardinality~$m$\rom;
\item\label{sets.dif}
if $r,s\in\fS$, then $\ls|r\sdif s|\in\{0, 4, 6, 8\}$
\rom($\sdif$ being the symmetric difference\rom).
\endroster
Then, for small values $(n,m)$, the maximum $\CA_{m,n}:=\max\ls|\fS|$
is as follows\rom:
\[*
\minitab\
(n,m):&    (n,1)&(n,2)             &(6,3)&(7,3)&(8,3)&(9,3)&(8,4)\cr
\CA_{m,n}:&    1&\lfloor n/2\rfloor&    4&    7&    8&   12&   14\cr
\endminitab
\]
More generally,}
\[*
{\4
\CA_{3,n}\le\left\lfloor\frac{n}3\Bigl\lfloor\frac{n-1}2\Bigr\rfloor\right\rfloor;\qquad
\CA_{m,n}\le\left\lfloor\frac1m\binom{n}{m - 1}\right\rfloor\quad\mbox{for $m\ge1$}.}
\]
\endlemma

\lemma[see~\cite{DIO}]\label{lem.max.set.D}
{\4The maximal cardinality
of a collection~$\fS$ of subsets $s\subset S$ satisfying condition~\iref{sets.dif}
of \autoref{lem.max.set.A}
is bounded \via\
\[*
\ls|\fS|\le\max_{m\ge0}
 \bigl(\CA_{m,n}+\CA_{m+2,n}+\CA_{m+4,n}+\CA_{m+6,n}+\CA_{m+8,n}\bigr),
\]
where
$\CA_{m,n}$ is as in \autoref{lem.max.set.A} and we let $\CA_{m,n}=0$
unless $0\le m\le n$.}
\endlemma

{\4The estimates on the number of integral vectors in
$\bA$--$\bD$--$\bE$ lattices also change due to~\eqref{eq.intr},
but they are quite straightforward.
Details will appear in~\cite{DIO}.
}

\subsection{Patterns\pdfstr{}{ \rm(see~\cite{degt:sextics})}}\label{s.patterns}
Let $\Cluster\subset\Orb$ be a
subset.
Then,
\roster*
\item
a \emph{pattern} is a function $\pat\:\Cluster\to\NN$ such that
$\pat(\orb)\in\sbnd(\orb)$ for each $\orb\in\Cluster$;
\item
a set $\fL\subset\fF$ is said to \emph{fit} a pattern~$\pat$ if
$\ls|\fL\cap\orb|=\pat(\orb)$ for each $\orb\in\Cluster$;
\item
the \emph{stabilizer} $\stab\pat$ is the stabilizer of~$\pat$ under
the action of
{\4the stabilizer $\stab\Cluster$ of~$\Cluster$ under $\stab\hbar$.}
\endroster
We use patterns to compute the sets $\BB_d(\Cluster)$.
We start with computing the
$(\stab\Cluster)$-orbits of the set of all patterns~$\pat$
satisfying the inequality $\sum\pat(\orb)\ge\bnd(\Cluster)-d$.
Then we pick a representative~$\pat$ and sort the combinatorial orbits
$\orb_1,\ldots,\orb_N\in\Cluster$ by the decreasing of the value of~$\pat$.
A geometric set~$\fL$ fitting~$\pat$ is constructed orbit-by-orbit,
$\varnothing=\fL_0\subset\fL_1\subset\ldots\subset\fL_N$. At each step~$k$,
we proceed as follows:
\roster*
\item
compute the $\RGL[\fL_{k-1}]$-orbits of sets $\fL'\in\geom(\orb_k)$ of size
$\ls|\fL'|=\pat(\orb_k)$;
\item
pick a representative $\fL'$ of each orbit and let
$\fL_k:=\sat(\fL_{k-1}\cup\fL')$;
\item
check that $\fL_k$ is geometric and that $\ls|\fL_k\cap\orb_i|=\pat(\orb_i)$
for $i\le k$.
\endroster
From this point on, the orbits $\orb_1,\ldots,\orb_k$ are considered
\emph{frozen}: for all subsequent extensions $\fL'\supset\fL_k$ we require that
$\fL'\cap\orb_i=\fL_k\cap\orb_i$, $1\le i\le k$. In other words, we decorate
$\fL_k$ with the \emph{frozen pattern}
$\anchor[@]{patL}\pat_{\fL_k}\:\{\orb_1,\ldots,\orb_k\}\to\NN$,
$\orb\mapsto\ls|\fL_k\cap\orb|$, which is to be enforced in all subsequent
extensions.

We run this algorithm in the hope that none of the chains survives through the
last step, yielding that $\BB_d(\Cluster)=\varnothing$. Otherwise,
the sets $\fL$ found are those generated by $\fL\cap\Cluster$;
they need to be examined
by other means, \cf. \autoref{s.extra.lines}, \autoref{s.max} below.

\remark\label{rem.comma}
Typically,
$\Cluster=\borb_1\cup\ldots\cup\borb_N$
is a union of orbits; they are processed one at a time,
and the computation (sorting of the combinatorial orbits)
depends on their order. For this reason,
we use comma separated lists: $\BB_d(\borb_1,\ldots,\borb_N)$.
On a few occasions, two orbits $\borb'$, $\borb''$ with
$\sbnd(\orb')=\sbnd(\orb'')$, $\orb^*\in\borb^*$, are processed
at once; this fact is indicated by $\borb'\orbitplus\borb''$ in the notation.
\endremark

\remark\label{s.check}
If $\bnd(\Orb)$ exceeds the goal~\eqref{eq.threshold}
by just a few units, we
use patterns to show directly that $\Bnd_{\bnd(\Orb)-M}(\Orb)=\varnothing$.
These cases are
marked with a \checkmark\ in the tables, and
any further explanation, as well as the list of orbits, are omitted.
\endremark

\subsection{Clusters and mega-clusters\pdfstr{}{ \rm(see~\cite{degt:sextics})}}\label{s.clusters}\label{s.mega}
If an orbit $\borb\subset\Orb$ is
large, we cannot effectively compute all patterns $\pat\:\borb\to\NN$. In
this case, we partition~$\borb$ into a number of \emph{clusters}~$\cluster_k$,
not necessarily disjoint, and proceed cluster-by-cluster.
Clusters are chosen to constitute a single $(\stab\hbar)$-orbit and,
when building a set~$\fL$, we
order them lexicographically, by the decreasing of the
{\4pair $\bigl(\ls|\fL\cap\cluster_k|,\nu(\cluster_k,\fL)\bigr)$, where}
\[*
{\4\nu(\fL,\cluster_k):=(\nu_0,\nu_1,\ldots),
\qquad
\nu_i:=\#\bigl\{\orb\in\cluster_k\bigm|\ddefect(\fL;\orb)=i\bigr\}.}
\]
Denoting by~$N$ the number of clusters and by~$m$, the \emph{multiplicity} of
the partition (the number of clusters containing each combinatorial orbit
$\orb\in\borb$), this convention implies that, for $\fL\in\BB_d(\borb)$, one
must have $\ddefect(\fL;\cluster_1)\le md/N$. More generally, at a
step~$(k+1)$,
\[*
(N-k)\ddefect(\fL;\cluster_{k+1})\le
 md-\sum_{i=1}^k\ddefect(\fL;\cluster_i).
\]
At this step,
one should try for~$\cluster_{k+1}$ a representative of
\emph{each $(\stab\patL\fL)$-orbit}.

\example\label{ex.clusters}
{\4As a purely artificial example, assume that an orbit $\borb$ of size $60$
is partitioned into six disjoint ($m=1$) clusters $\cluster_k$, each consisting of
$10$ combinatorial orbits, and we are trying to construct a set
$\fL\in\BB_{15}(\borb)$. (A direct computation of all patterns on~$\borb$ is
clearly not feasible.) By the $(\stab\hbar)$-symmetry, we can assume that
$\ddefect(\fL,\cluster_1)\le2$ and easily compute
$\BB_2(\cluster_1)$. For each set $\fL_1\in\BB_2(\cluster_1)$, we
compute the orbits of the action of $\stab\patL{\fL_1}$ on the remaining five
clusters. At this point, the process bifurcates: for \emph{each} orbit,
we pick a representative, denoted by $\cluster_2$, and run
a separate branch of the algorithm with this representative. An extension
$\fL_2\supset\fL_1$ is constructed using a pattern $\pat_2\:\cluster_2\to\NN$
subject to the following conditions (recall that larger clusters are to be
added first):
\roster*
\item
if $\ddefect(\fL,\cluster_1)=0$, then $0\le\ddefect(\fL,\cluster_2)\le3$;
\item
if $\ddefect(\fL,\cluster_1)=1$, then $1\le\ddefect(\fL,\cluster_2)\le2$;
\item
if $\ddefect(\fL,\cluster_1)=2$, then $\ddefect(\fL,\cluster_2)=2$.
\endroster
We continue with~$\cluster_3$, $\cluster_4$, \etc. in a similar manner;
usually, this process terminates
long before all (six) clusters have been used, see \autoref{s.extra.lines}
below.

The order of clusters is used to further reduce the overcounting.
For example, if $\nu(\fL,\cluster_1)=(4,2,\ldots)$, then
$\nu(\fL,\cluster_2)\ne(5,0,1,\ldots)$: otherwise, $\cluster_2$ must be added
first. This observation reduces the number of patterns that can be used to
build~$\fL_2$.}
\endexample

If the number of clusters is still large, our algorithm tends to diverge due
to the massive repetition caused by the random choice of subsequent clusters.
To remedy this, we partition the set of clusters
into \emph{mega-clusters}~$\mega_n$, which are also assumed to constitute a
single $(\stab\hbar)$-orbit. Arguing as above, we fill, cluster-by-cluster,
the \emph{first} mega-cluster~$\mega_1$, \emph{assumed maximal}. Usually, the
algorithm stops at this point; in a few cases, one needs to extend the sets
obtained by one extra
\emph{cluster}.

\example\label{ex.mega}
{\4Assume that an orbit $\borb$ is partitioned into $60$
\emph{clusters} $\cluster_k$ which are grouped into six disjoint
\emph{mega-clusters} $\mega_n$, and we are trying to construct a set
$\fL\in\BB_{15}(\borb)$.
The computation runs exactly as in \autoref{ex.clusters},
but
instead of building $\fL\cap\cluster_k$ from the precomputed sets
$\fL\cap\orb_i$, we build $\fL\cap\mega_n$ from $\fL\cap\cluster_k$ (which, in
turn, are constructed on the fly from $\fL\cap\orb_i$).
In other words,
instead of starting from at least $45$ maximal (hence,
indistinguishable by their patterns) clusters picked, essentially, in a
random order, we assert that $\ddefect(\fL,\mega_1)\le2$ and build the
(almost maximal) set $\fL\cap\mega_1$. As stated, in the applications
the process tends to
terminate at this point, due to \autoref{s.extra.lines} below.
}
\endexample

\subsection{Extension by extra conics}\label{s.extra.lines}
According to \autoref{prop.I}\iref{I.rank} (and to our choice to consider
primitive sublattices $S\subset\bN$ only), any saturated set $\fL\subset\fF$ of
rank~$20$ is \emph{maximal}
in the sense that it has no proper geometric extensions.
Hence, if such a set $\fL:=\fL_k$ appears, at any step, in the algorithm
of \autoref{s.patterns}, it can be \emph{discarded} immediately:
{\4that is, we check if $\fL$ satisfies~\eqref{eq.threshold} (and record it as
an exception if it does), but
\emph{we do not continue the current algorithm on~$\fL$}, \ie, we do not try
to construct further extensions $\ldots\supset\fL_{k+1}\supset\fL_k$.}
(In {\4practice}, in order to obtain plenty of
examples, we recorded all geometric sets $\fL$ satisfying $\ls|\fL|>200$.
{\4However, we do not assert that this list is complete.})

To save time, we push this policy two steps further and discard,
\emph{as soon as they appear},
all sets~$\fL$ of rank $\rank\fL\ge18$. Thus, we compute
and use in \autoref{lem.bounds}
the sets
\[*
\anchor[@]{B''}
\BBii_d(\Cluster):=\bigl\{\fL\in\BB_d(\Cluster)\bigm|\rank\fL<18\bigr\}.
\]
{\4More precisely, as soon as a set $\fL:=\fL_k$ of rank $\rank\fL=19$ has
been obtained at an intermediate step of the current algorithm, we
check~\eqref{eq.threshold} for $\fL$ itself and, instead of continuing the
algorithm, merely consider all corank~$1$ extensions $\fL'\supset\fL$
obtained by adding to~$\fL$ an extra conic. To this end,}
we partition the set $\fF\sminus\fL$ by the equivalence
relation
\[*
l'\sim_\fL l''\quad\text{iff}\quad\spn(\fL\cup l')=\spn(\fL\cup l'').
\]
Let $\Classesi(m)$,
$m\in\NN$, be the set of equivalence classes~$\class$ such that
\roster
\item\label{ext.root}
$\spn(\fL\cup\class)$ is root free
(otherwise, $\fL\cup\class$ is not admissible),
\item\label{ext.frozen}
$\class$ is disjoint from all frozen orbits
(so that $\fL\cup\class$ fits $\patL\fL$), and
\item\label{ext.size}
$\ls|\fL|+\ls|\class|\ge m$.
(We need $m=M$, but we use a weaker bound $m=200$.)
\endroster
Note that, since $\spn({-})$ is a primitive sublattice, this
computation can be done over~$\Q$ rather than~$\Z$, which makes it much
faster.
Then, we pick a representative~$\class$
of each $\OGL$-orbit of $\Classesi(m)$
and record $\fL\cup\class$ as an exception if it is geometric.

A similar approach is used if $\rank\fL=18$: here, in addition to corank~$1$,
we also check corank~$2$ extensions of~$\fL$. Technically,
we consider the unions $\class_{ij}:=\class_i\cup\class_j$ over all pairs
$\class_i\ne\class_j\in\Classesi(0)$ and proceed as follows:
\roster*
\item
declare $\class_{ij}\sim_\fL\class_{kl}$ whenever
$\spn(\fL\cup\class_{ij})=\spn(\fL\cup\class_{kl})$;
\item
replace each $\sim_\fL$-equivalence class with its \emph{union} and consider the
set~$\Classesii$ of these unions (the elements of~$\Classesii$ are
\emph{not} disjoint subsets of $\fL'\supset\fL$);
\item
let $\Classesii(m):=\{\class\in\Classesii\,|\,
\text{$c$ satisfies the obvious analogues of \iref{ext.root}--\iref{ext.size}}\}$;
\item
pick a representative~$\class$ of each $\OGL$-orbit
of~$\Classesii(m)$
and record $\fL\cup\class$ as an exception if it is geometric.
(Again, we use $m=200$.)
\endroster

\remark\label{rem.extension}
{\4This extension procedure gives us all sets
$\fL\in\BB_d(\Cluster)\sminus\BBii_d(\Cluster)$ or geometric oversets thereof
that satisfy~\eqref{eq.threshold}. That is why, when using
\autoref{lem.bounds}, we can replace the hypotheses
$\BB_d(\Cluster)=\varnothing$ with $\BBii_d(\Cluster)=\varnothing$.
The same observation applies to the \latin{a priori} even smaller sets
$\BBs_d(\Cluster)$ considered in \autoref{s.max} below.}
\endremark

\subsection{Extension by an extra orbit}\label{s.max}
In \autoref{s.extra.lines}, sets $\fL$ of rank $\rank\fL\ge18$ are discarded
\emph{on the way}, as soon as they appear in the course of the computation.
This results in the set $\BBii_d(\Cluster)$, which may still be non-empty.
(We tweak the parameters so that these sets are not very large, up to a dozen
of elements.)
{\4To find all geometric oversets of a set $\fL\in\BBii_d(\Cluster)$ that
satisfy~\eqref{eq.threshold} (\cf. \autoref{rem.extension}), we use the}
\emph{extension by a maximal orbit}.
Denote by $\bar\Orb_\fL$ the set of all combinatorial orbits~$\orb$
that are \emph{not}
frozen in $\fL$ and such that $\ls|\fL\cap\orb|<\bnd(\orb)$.
We assume that
\[*
\ls|\fL|+\sum_{\orb\in\bar\Orb_\fL}\bigl(\bnd(\orb)-\bnd_2(\orb)\bigr)<M,
\]
where $\bnd_2(\orb)$ is the second maximal element of $\sbnd(\orb)$.
({\4Below,} in the
only case where this assumption does not hold, we use the
\emph{extension by an arbitrary orbit}, which is quite similar.)
The inequality means that, for each geometric extension $\fL'\supset\fL$
satisfying~\eqref{eq.threshold},
at least one orbit
$\orb\in\bar\Orb_\fL$ must be maximal, $\ls|\fL'\cap\orb|=\bnd(\orb)$.
Thus, we can pick a representative $\orb$ of each orbit of the
$(\stab\patL\fL)$-action on $\bar\Orb_\fL$, extend
the frozen pattern~$\patL\fL$ by a single
value $\orb\mapsto\bnd(\orb)$, and run one more step of the algorithm of
\autoref{s.patterns}, followed by discarding,
{\4as explained in \autoref{s.extra.lines},}
all newly discovered sets
$\fL'\supset\fL$
of rank $\rank\fL'\ge18$.

Since we add at least one extra conic $l\notin\fL$, the rank of the set must
increase.
It follows that this procedure produces all geometric extensions of
a set $\fL\in\BBii_d(\Cluster)$ of rank $\rank\fL=17$,
and we apply it automatically, indicating this fact in the notation
\[*
\anchor[@]{B*}
\BBs_d(\Cluster):=\bigl\{\fL\in\BB_d(\Cluster)\bigm|\rank\fL<17\bigr\}.
\]
(In particular, the reference to this notation implies that
$\BBii_d(\Cluster)\ne\varnothing$.) In the very few cases where
$\BBs_d(\Cluster)\ne\varnothing$, the procedure needs to be applied several
times; we mention each of these {\4exceptional} cases explicitly.

\subsection{Single orbits\pdfstr{}{ \rm(see~\cite{DIO})}}\label{s.single}
We
denote by $\ssingle\subset\Orb$ the set of all combinatorial orbits
consisting of a single conic.
The sets of the form $\Bnd_d(\frak{s})$, $\frak{s}\subset\ssingle$,
are computed ``backwards'', \via\ \emph{iterated index~$2$ subgroups} (unless
$\ls|\single|<16$, in which case we merely analyze the $(\stab\hbar)$-orbits
of all combinations of~$\single$).

{\4Let $\ts:=\sat_\Z\frak{s}$ and $\spnc\fL:=\spn\fL\cap\spn_\Z\ts$ for
$\fL\subset\ts$.
Given}
a $\Z$-saturated subset $\fL\subset\4\ts$, we {\4can} compute the
orbits of the $\OGL$-action on the $\F_2$-\penalty0vector space
$(\spn_\Z\fL)\dual\otimes\F_2$ and, for a representative~$v$ of each orbit,
consider the
{\4annihilator}
$\fL_v:=\{l\in\fL\,|\,v(l)=0\}$. (The sets $\fL_v$ are
$\Z$-saturated but, in general, not saturated; still, they are retained at
the intermediate steps.)
{\4It is immediate that, starting from $\fL=\4\ts$ and iterating this
procedure, we obtain all \emph{relatively saturated} subsets
$\fL\subset\4\ts$
(\ie, such that
$\fL=\fF\cap\spnc\fL$),
and there remains to
select those of size $\ls|\fL|\ge\ls|\frak{s}|-d$ and whose saturation is
geometric.
Indeed, if $\fL$ is relatively saturated and $\fL\ne\ts$, then
$\spn_\Z\fL\otimes\F_2\subset\spn_Z\ts\otimes\F_2$ is a proper subspace and
there is a covector $v\ne0$ vanishing on~$\fL$. Then, $\ts_v\supset\fL$ is a
proper subset of~$\ts$, and we can replace $\ts$ with $\ts_v$ and argue by
induction on the size $\ls|\ts\sminus\fL|<\infty$.}



\remark\label{rem.single}
When computing a set of the form $\BBii_d(\single,\borb_1,\ldots)$, we always
start with $\BB_d(\single)$, using iterated index~$2$ subgroups; then,
all combinatorial orbits $\orb\in\single$ are considered frozen, and we
continue with the
algorithm of \autoref{s.patterns}.
\endremark

\subsection{Replanting}\label{s.replant}
Let $\bN\ni\hbar$ be a $12$-polarized Niemeier lattice.
Let $\fF:=\fF(\bN,\hbar)$, denote
$F:=\spn_\Z\fF$, and assume that $\rank F=24$, so that
$F$ is of finite index in~$\bN$. Since $\bN$ is unimodular and
$\hbar\in4F\dual$, we have $[\bN:F]\ge4$.
Till the end of this section,
{\em we assume that $[\bN:F]=4$ and, hence, $\ls|\discr F|=16$}.

Since $\bN$ is unimodular and $\hbar\in\bN$ is a primitive vector, there is a
vector $a\in\bN$ such that $a\cdot\hbar=1$, and it is immediate that
$\discr F\cong\Z/4\oplus\Z/4$ is generated by the order~$4$ classes
$\Ga:=a\bmod F$ and $\eta:=\frac14\hbar$.
This group has but two isotropic subgroups of order~$4$, both cyclic,
generated by $\Ga$ or $\Ga':=\Ga+2\eta$.
Hence, $F$ has two unimodular finite index extensions: $\bN$ (defined
by~$\Ga$) and another $12$-polarized
Niemeier lattice $\bN'\ni\hbar$ defined by~$\Ga'$.
We say that $\fF(\bN,\hbar)$ is \emph{replanted} to $\bN'$.
The next lemma states that, for the proof of \autoref{th.main}, it suffices
to consider only one of the two lattices.

\lemma\label{lem.replant}
There is a canonical bijection between the admissible/geometric sets
in $\bN\ni\hbar$ and those in $\bN'\ni\hbar$.
\endlemma

\proof
For any $m\in\Z/4\sminus0$, we have $(m\Ga)\cdot\hbar\ne0\bmod4$ and
$(m\Ga')\cdot\hbar\ne0\bmod4$.
Hence, $\fF(\bN',\hbar)=\fF$, \ie, the two lattices share the same set of
{\4conics}. Furthermore,
\emph{if $\fL\subset\fF$ is admissible in~$\bN$} (or~$\bN'$), then
necessarily $\spn\fL\subset F$ (as one must have $\hbar\in4(\spn\fL)\dual$),
\ie, $\fL$ has the same span in both lattices.
All other conditions are stated in terms of $\spn\fL$.
\endproof

There also are canonical isomorphisms $\OG_\hbar(\bN)=\OG_\hbar(\bN')$
and $\RG_\hbar(\bN)=\RG_\hbar(\bN')$ (for the latter, any root
\emph{orthogonal to~$\hbar$} must lie in~$F$); hence, the two lattices have
the same structure of orbits and combinatorial orbits.
Note also that $\bN$ and $\bN'$ share a common index~$2$ extension of~$F$, \viz.
the one defined by $2\Ga=2\Ga'$.

\makeatletter
The pair $\bN\ni\hbar$ is called \emph{reflexive} if
$(\bN',\hbar)\cong(\bN,\hbar)$ (where we still assume that $[\bN:F]=4$).
If this is the case,
any isomorphism $\Gf\:(\bN',\hbar)\to(\bN,\hbar)$ restricts to an
automorphism of $F\ni\hbar$ that does not extend to~$\bN$.
The induced permutation~$\anchor[@]{tau}\pperm$
of the set of $\OG_\hbar(\bN)$-orbits
of~$\fF$ is an involution
independent of the choice of~$\Gf$; it is called the
\emph{replanting involution}.
\makeatother

\section{Lattices with few components}\label{S.few}

Recall (see~\cite{Niemeier} or \cite{Conway.Sloane}) that $23$ of the $24$
Niemeier lattices are rationally generated by roots.
We use the notation
\roster*
\item
$\bN:=\N(\bR)$, where $\bR\subset\bN$ is the maximal root system,
\item
$\bR=\bigoplus_{k\in\Omega}\bR_k$ is the decomposition
into irreducible components, and
\item
we reserve the
notation~$\Omega$ for the index set.
\endroster
For a vector $v\in\bN$, we write $v=\sum_kv_k$, where
$v_k\in\bR_k\dual$, $k\in\Omega$, is the orthogonal projection.
The \emph{support} of~$v$ is $\supp v:=\{k\in\Omega\,|\,v_k\ne0\}$.


In this section, we treat the $21$ Niemeier lattices with $\ls|\Omega|\le8$.

\theorem\label{th.<=8}
Let $\bN$ be a Niemeier lattice generated over~$\Q$ by a root system~$\bR$ with
at most eight irreducible components. Then, for any square~$12$ vector
$\hbar\in\bN$ and any geometric set $\fL\in\fF(\bN,\hbar)$ one has
$\ls|\fL|\le260$.
\endtheorem

\proof
For each lattice $\bN:=\N(\bR)$, we list the $\OG(\bN)$-orbits of square~$12$
vectors $\hbar\in\bN$, pick a representative of each orbit, and
consider pairs $\bN\ni\hbar$ one-by-one.

For each pair $\bN\ni\hbar$, we list orbits~$\borb_n$ and combinatorial
orbits~$\orb$ and, for the latter, use~\cite{DIO}
{\4(see also \autoref{s.estimates})}
to compute rough
combinatorial estimates $\bnd'(\orb)\ge\bnd(\orb)$. In most cases, we have
$\bnd(\Orb)\le\bnd'(\Orb):=\sum_nm(\borb_n)\bnd'(\orb)<M$, $\orb\in\borb_n$,
where
$m(\borb_n)$ is the number of combinatorial orbits $\orb\in\borb_n$.
If not, we compute the sharp bounds $\bnd(\orb)$ (by brute force; in
general, they are less than those in~\cite{DIO} as we consider primitive
sublattices of~$\bN$ only). If $\bnd(\Orb)<M$, we disregard the pair
$\bN\ni\hbar$. In particular, for $19$ lattices we have $\bnd(\Orb)\le258$
for any vector $\hbar\in\bN$.

The few remaining pairs are listed in the tables below. For each
vector~$\hbar\in\bN$, we show the numbers $\ls|\fF|$ and
$\bnd'(\Orb)\to\bnd(\Orb)$ and list the
orbits $\borb_n$
(only those used in the computation),
indicating for each the multiplicity $m(\borb_n)$ and
count $\ls|\orb|$ and bounds $\bnd'(\orb)\to\bnd(\orb)$,
$\orb\in\borb_n$. (The bounds known to be sharp are underlined.)

\convention\label{conv.tables}
For the components $\hbar_k$ of~$\hbar$ we use the notation
$[\![\hbar_k^2]\!]_d$,
where $d$ is either the discriminant class
of~$\hbar_k^2$ (in the notation of \cite{Conway.Sloane}) or,
if $\hbar_k\in\bR_k$, the symbol
\[*
\mbox{$0$ (if $\hbar_k=0$),\quad
$\circ$ (if $\hbar_k^2=2$),\quad
$\bullet$ (if $\hbar_k^2=4$),\quad
$*$ (if $\hbar_k^2=6$)}.
\]
If these data do not determine~$\hbar_k$, we use an extra superscript
whose precise meaning is not very important (see~\cite{degt:sextics}).

For the components $l_k$ of a conic, we use the notation
$[l_k\cdot\hbar_k]_d$, where $d$ and an occasional superscript have
the same meaning as for~$\hbar$.
(It is worth emphasizing that each of $d:=d_k$, $l_k^2$, and
$l_k\cdot\hbar_k$, $k\in\Omega$, is constant within each combinatorial
orbit~$\orb$; in the tables, we show one representative $\orb\in\borb_n$
only.)


When describing partitions, in addition to $\supp\orb:=\supp l$, $l\in\orb$,
which is also constant,
we use the notation
\[*
\anchor[@]{pos}
\ppos(q):=\bigl\{k\in\Omega\bigm|\hbar_k^2=q\bigr\},\quad q\in\Q,
\]
and $\pos_m(q)$ for the set of $m$-combinations of $\pos(q)$, $m\ge2$.
\endconvention


\removelastskip
\input{\listdirectory/\prefix_D24}
\input{\listdirectory/\prefix_D16_E8}
\input{\listdirectory/\prefix_3E8}
\input{\listdirectory/\prefix_A24}
\input{\listdirectory/\prefix_2D12}
\input{\listdirectory/\prefix_A17_E7}
\input{\listdirectory/\prefix_D10_2E7}
\input{\listdirectory/\prefix_A15_D9}
\input{\listdirectory/\prefix_3D8}
\input{\listdirectory/\prefix_2A12}
\input{\listdirectory/\prefix_A11_D7_E6}
\input{\listdirectory/\prefix_4E6}
\input{\listdirectory/\prefix_2A9_D6}
\input{\listdirectory/\prefix_4D6}
\input{\listdirectory/\prefix_3A8}
\input{\listdirectory/\prefix_2A7_2D5}
\input{\listdirectory/\prefix_4A6}
\input{\listdirectory/\prefix_4A5_D4}

\checkline[6D4]{1}{\hidelines}
\input{\listdirectory/\prefix_6D4}
The only configuration with
$\bnd(\Orb)\ge261$ is replanted (see \autoref{s.replant})
to \config[24A1]{11} in $\N(24\bA_1)$ and is considered
in \autoref{s.24A1-11} below.

\input{\listdirectory/\prefix_6A4}

\setlattice{8A3}
\checkline{1}{\hidelines}
\checkline{2}{\hidelines}
\checkline{4}{\checkdone}
\checkline{9}{\hidelines}
\checkline{10}{\hidelines}
\checkline{13}{\hidelines}
\checkline{17}{\checkdone}
\checkline{19}{\hidelines}
\checkline{21}{\hidelines}
\checkline{24}{\hidelines}
\checkline{26}{\checkdone}

\input{\listdirectory/\prefix_8A3}
Ten of the eleven relevant configurations
(those with $\bnd(\Orb)\ge261$) are either replanted to other
lattices or covered by \autoref{s.check}.

\subsubsection{Replanted configurations}
The following seven configurations are replanted (see \autoref{s.replant}) to other
lattices and, hence, considered elsewhere:
\roster*
\hitem{1} is replanted to
\config[24A1]{7} in $\N(24\bA_1)$, see \autoref{s.24A1-7} below;
\gchecked
\hitem{2} is replanted to
\config[24A1]{12} in $\N(24\bA_1)$, see \autoref{s.24A1-12} below;
\gchecked
\hitem{9} is replanted to \config[12A2]{14} in $\N(12\bA_2)$,
see \autoref{s.12A2-14} below;
\gchecked
\hitem{10} is replanted to \config[12A2]{15} in $\N(12\bA_2)$,
see \autoref{s.12A2-15} below;
\gchecked
\hitem{13} is replanted to \config[12A2]{16} in $\N(12\bA_2)$,
see \autoref{s.12A2-16} below;
\gchecked
\hitem{19} is replanted to \config[12A2]{17} in $\N(12\bA_2)$,
see \autoref{s.12A2-17} below;
\gchecked
\hitem{21} is replanted to \config[12A2]{18} in $\N(12\bA_2)$,
see \autoref{s.12A2-18} below.
\gchecked
\endroster

\def\porb#1{\obj5{\cluster_{#1}}}
\subsubsection{Vector~\hlink{8A3}{24}\rm:}\label{s.8A3-24}
we have
$\Bnd_5(\single)=
\Bii11[ 5 ]=
\varnothing$;
for the latter, we use the partition
\[*
\porb{k}:=\bigl\{\orb\subset\oorb5\bigm|k\notin\supp\orb\bigr\},
\quad k\in\supp\hh.
\gchecked
\]
This case completes the proof of \autoref{th.<=8}.
\endproof

\setlattice{12A2}

\section{The lattice $\N(12\bA\sb2)$}\label{S.12A2}

The goal of this section is the following theorem.

\theorem\label{th.12A2}
Let $\bN:=\N(12\bA_2)$ and $\hbar\in\bN$, $\hbar^2=12$.
Then, with one exception \rom(up to automorphism\rom),
for any geometric set $\fL\subset\fF(\bN,\hbar)$, one has $\ls|\fL|\le260$.
The exception is the set $\Lmax1$ of size~$285$, see~\eqref{eq.Lmax.1} below.
\endtheorem

\checkline{1}{\hidelines}
\checkline{3}{\hidelines}
\checkline{4}{\hidelines}
\checkline{5}{\hidelines}
\checkline{6}{\hidelines}
\checkline{7}{\hidelines}
\checkline{8}{\hidelines}
\checkline{9}{\hidelines}
\checkline{10}{\hidelines}
\checkline{11}{\hidelines}
\checkline{12}{\hidelines}
\checkline{13}{\hidelines}
\checkline{14}{\hidelines}
\checkline{15}{\hidelines}
\checkline{16}{\hidelines}
\checkline{17}{\hidelines}
\checkline{18}{\hidelines}

\breaktable{9}
\breaktable{13}


\proof\leavevmode
{\def\vectorbox#1{#1\,\,}\def\countskip{\ }\def\ltitle{}\input{\listdirectory/\prefix_12A2}}
{\4As above, we list
only those vectors~$\hbar$ for which
$\bnd(\Orb)\ge M$, \cf. also \autoref{s.check}.}

\subsection{Replanted configurations}
Two configurations can be replanted (see \autoref{s.replant}) to $\N(24\bA_1)$
and are considered in \autoref{S.24A1}. They are as follows:
\roster*
\hitem{1}
is replanted to \config[24A1]{8} in $\N(24\bA_1)$,
see \autoref{s.24A1-8} below;
\gchecked
\hitem{3} is replanted to \config[24A1]{9} in $\N(24\bA_1)$,
see \autoref{s.24A1-9} below.
\gchecked
\endroster

\subsection{Vector~\hlink{12A2}{2}}\label{s.12A2-2}
We have $\ls|\stab\hh|=720$ and $\rank\fF=21$.
Since the first combinatorial estimate
$\bnd(\orb)\le27$, $\orb\in\oorb1$, is too high, we
use \emph{blocks} (see \cite[\S\,4.2]{degt:sextics} for details) to
compute \emph{large} geometric intersections $\fL\cap\orb$ only, arriving at
\[*
\text{$\ls|\fL\cap\orb|=19$ or $\ls|\fL\cap\orb|\le15$
 for any $\orb\in\oorb1$ and geometric set $\fL\subset\fF$}.
\]
It follows that $\ls|\fL|\le245$ whenever
$\ls|\fL\cap\orb|<19$ for at least ten
orbits $\orb\in\oorb1$.
Then, using patterns (see \autoref{s.patterns}), we compute the geometric
sets~$\fL$ satisfying $\ls|\fL\cap\orb|=19$ for (at least) six
orbits $\orb\in\oorb1$. We find one set $\Lmax1$ of rank~$20$ and size~$285$
and one set~$\fL$ of rank~$18$; extending the latter as in
\autoref{s.extra.lines}, we obtain, apart from $\Lmax1$, five sets of sizes
$249$, $231$, $213$, $213$, $207$.

The maximizing set $\Lmax1$ can be described as
\[
\Lmax1=\fF\cap r^\perp,\quad
r:=\sum_{k\in\Omega}(\text{a root $r_k\in\hh_k^\perp\subset\bR_k$}).
\maxlabel1
\gchecked
\]
Clearly, vectors $r_k\in\bR_k$ and $r$ as
in~\eqref{eq.Lmax.1} are unique up to the action
of $\RG_\hh$.

Another large set worth mentioning is
\[
\Lmisc1=\fF\cap(r+r_k)^\perp,\quad k\in\Omega\sminus\supp\hh,
\misclabel1{249}
\]
where $r_k\in\bR_k$ and $r=\sum r_k$ are as in~\eqref{eq.Lmax.1}. This set is
discussed in \autoref{ex.249}.


\def\porb#1{\obj{2}{\cluster_{#1}}}
\def\pporb#1{\obj{5}{\cluster'_{#1}}}
\def\ppporb#1{\obj{6}{\cluster''_{#1}}}
\subsection{Vector~\hlink{12A2}{4}}\label{s.12A2-4}
Using the partitions
\begin{alignat*}3
\porb{k}&:=\bigl\{\orb\in\oorb2\bigm|l_k\cdot\hh_k=-\tfrac13\bigr\},\quad
 &k&\in\pos\bigl(\tfrac23\bigr)&\quad&\text{($1$-fold)},\\
\pporb{k}&:=\bigl\{\orb\in\oorb5\bigm|l_k\cdot\hh_k=\tfrac23\bigr\},\quad
 &k&\in\pos\bigl(\tfrac83\bigr)&&\text{($1$-fold)},\\
\ppporb{k}&:=\bigl\{\orb\in\oorb6\bigm|
 l_k\cdot\hh_k=\tfrac13\bigr\},\quad
 &k&\in\pos\bigl(\tfrac23\bigr)
 &&\text{($2$-fold)},
\end{alignat*}
we show that
$\Bs24[0,2]=    
\Bii19[5]=
\Bii33[6]=
\varnothing$.
\gchecked

\def\porb#1{\obj{3}{\cluster_{#1}}}
\def\pporb#1{\obj{5}{\cluster'_{#1}}}
\subsection{Vector~\hlink{12A2}{5}}\label{s.12A2-5}
Using the partitions
\begin{alignat*}3
\porb{k}&:=\bigl\{\orb\in\oorb3\bigm|l_k\cdot\hh_k=\tfrac23\bigr\},\quad
 &k&\in\pos\bigl(\tfrac23\bigr)&\quad&\text{($2$-fold)},\\
\pporb{k}&:=\bigl\{\orb\in\oorb5\bigm|l_k\cdot\hh_k=-\tfrac13\bigr\},\quad
 &k&\in\pos\bigl(\tfrac23\bigr)&&\text{($1$-fold)},
\end{alignat*}
we show that $\Bii35[3]=\varnothing$ and compute
$\Bs18[0,5]$; the latter consists of four sets of size
$111$, $99$, $99$, $95$ and rank $16$, $16$, $15$, $16$, respectively,
which are extended
by a maximal orbit (see \autoref{s.max}) \emph{twice}.
Similarly, the only set $\fL\in\Bs17[7, 6, 1]$ has rank $16$ and can be
extended twice by a maximal orbit.
The largest set observed in this computation has size $249$;
letting $\{n\}:=\pos(2)$ and $\{p,q\}:=\pos(0)$, this set can be described
as
\[
\Lmisc2:=\fF\cap\spn\bigl(\hh-2\hh_n,r+r_p,r_n,r_q,\bigr)^\perp,
\misclabel2{249}
\]
where a collection of roots $r_k\in\hh_k^\perp\subset\bR_k$, $k\in\Omega$,
and vector $r:=\sum_{k\in\Omega}r_k$ are as in~\eqref{eq.Lmax.1}.
This set is discussed in \autoref{ex.249} below.
\gchecked

\subsection{Vector~\hlink{12A2}{6}}\label{s.12A2-6}
We have
\[*
\Bs12[8,16]=   
\Bs12[17+14]=  
\Bii12[0,6]=
\Bii9[12]=
\Bii3[7]=
\varnothing.
\]
To complete the
hypotheses of \autoref{lem.bounds},
we compute the set $\Bii7[2,11,-9]=\varnothing$
employing the replanting involution $\perm=(\oind2,9)$ (see
\autoref{s.replant}):
due to the symmetry,
we can assume that $\defect[\fL;2]\le\defect[\fL;-9]$;
hence, it suffices to take $\defect[\fL;2]\le3$ at step~$1$ and
$2\defect[\fL;2]+\defect[\fL;11]\le7$ at step~$2$.
\gchecked

\def\porb#1{\obj{7}{\cluster_{#1}}}
\subsection{Vector~\hlink{12A2}{8}}\label{s.12A2-8}
Using the partition
\[*
\porb{k}:=\bigl\{\orb\in\oorb7\bigm|l_k\cdot\hh_k=\tfrac23\bigr\},\quad
 k\in\pos\bigl(\tfrac23\bigr),
\]
we show that $\Bii12[7]=\varnothing$.
Besides, $\Bii11[6]=\varnothing$.
\gchecked

\subsection{Vector~\hlink{12A2}{9}}\label{s.12A2-9}
This set is reflexive (see \autoref{s.replant}), the replanting
involution being
\[*
\perm=(\oind{1},20)(2,19)(\oind{4},13)(7,\oind{17})(8,\oind{16}).
\]
We have
$\Bii10[6,11]=\Bii15[0,9]=\varnothing$ and
$\Bii5[4]=\Bs7[ 1, 17, 16 ]=\varnothing$;
from the latter, applying~$\perm$, we derive that also
$\Bii5[-13]=\Bs7[-20,-7,-8]=\varnothing$.
\gchecked

\subsection{Vector~\hlink{12A2}{10}}\label{s.12A2-10}
We have
$
\Bii5[5,12]=
\Bs8[3,6]=  
\Bs8[7,9]=  
\varnothing
$.
Then, using the partition
\[*
\porb{k}:=\bigl\{\orb\in\oorb8\bigm|k\notin\supp\orb\bigr\},\quad
k\in\pos(2),
\]
we show that $\Bii15[8]=\varnothing$ and
$\defect[\fL;8]>34$ for each $\fL\in\Bii18[0]$.
\gchecked

\def\porb#1{\obj{2}{\cluster_{#1}}}
\subsection{Vector~\hlink{12A2}{11}}\label{s.12A2-11}
We have
$
\B  17 [ 0 ]=
\Bii48 [ 2 ]=
\varnothing
$;
for $\oorb2$, we use the partition
\[*
\porb{c}:=\bigl\{\orb\in\oorb{2}\bigm|c\subset\supp\orb\bigr\},\quad
c\in\pos_2(0).
\gchecked
\]

\subsection{Other configurations}\label{s.12A2-other}
For the remaining eight vectors $\hbar\in\bN$, we directly use
patterns to compute a few sets of the form $\BBii_d(\Cluster)$ and refer to
\autoref{lem.bounds}.
\bgroup\let\linkform\configform\let\hitem\hitemii
\roster
\hitem{7}
$\Bii11[5]=
\Bii8[6]=
\varnothing$;
\gchecked

\hitem{12}
$\Bii15[5]=
\Bii9[6]=
\Bii5[8]=
\varnothing$;
\gchecked

\hitem{13}
$
\Bii13[10]=
\Bii8[8+11]=
\Bs9[0,13]= 
\varnothing
$;
\gchecked

\hitem{14}
$
\Bii5[18, 22,  4]=
\Bii5[ 8, 25, 37]=
\Bii5[23, 32]=
\Bs4[ 0, 21, 9 ]=   
\varnothing
$;
\gchecked

\hitem{15}
$
\Bii9[ 7 ]=
\Bii4[ 8 ]=
\varnothing
$;
\gchecked

\hitem{16}
$
\Bii5[ 10, 11 ]=
\Bii4[ 19, 12 ]=
\Bii4[ 15, 8 ]=
\varnothing
$;
\gchecked

\hitem{17}
$
\Bii6[ 14, 24 ]=
\Bs12[ 0, 17, 33, 9 ]=  
\varnothing
$;
\gchecked

\hitem{18}
$
\Bii5[ 14, 20 ]=
\Bii5[ 0, 17 ]=
\Bs 7[ 9 ]=     
\varnothing
$.
\gchecked

\endroster
This completes the proof of \autoref{th.12A2}.
\egroup
\endproof

\makeatletter
\let\@int\@@int
\let\@rat\@@rat
\makeatother

\setlattice{24A1}
\checkline{2}{\hidelines}
\checkline{3}{\hidelines}
\checkline{4}{\hidelines}
\checkline{5}{\hidelines}
\checkline{6}{\hidelines}
\checkline{7}{\hidelines}
\checkline{8}{\hidelines}
\checkline{9}{\hidelines}
\checkline{11}{\hidelines}
\checkline{12}{\hidelines}
\def\ltitle{}

\section{The lattice $\N(24\bA\sb1)$}\label{S.24A1}

The goal of this section is the following theorem.

\theorem\label{th.24A1}
Let $\bN:=\N(24\bA_1)$ and $\hbar\in\bN$, $\hbar^2=12$.
Then, with few exceptions, one has $\ls|\fL|\le260$
for any geometric set $\fL\subset\fF(\bN,\hbar)$.
Up to automorphism, there are six exceptional sets\rom:
\roster*
\item
three sets
$\Lmax2$, $\Lmax3$, $\Lmax4$ of size~$285$,
see \eqref{eq.Lmax.2}, \eqref{eq.Lmax.3}, \eqref{eq.Lmax.4} below\rom;
\item
three sets
$\Lsub1$, $\Lsub2$, $\Lsub3$ of size~$261$,
see \eqref{eq.Lsub.1}, \eqref{eq.Lsub.2}, \eqref{eq.Lsub.3} below.
\endroster
\endtheorem

\proof
\input{\listdirectory/\prefix_24A1}
{\4(As above, we list
only those vectors~$\hbar$ for which
$\bnd(\Orb)\ge M$, \cf. also \autoref{s.check}.)}
To save space, we fix a root basis $\{r_k\}$, $k\in\Omega$, for
$\bR=24\bA_1$
and, given a vector $v\in\bN$ (either $\hbar$ or a conic~$l$),
use the following notation for its components $v_k=\Ga r_k\in\bR_k\dual$:
\[*
\mbox{$\Azero$ ($\Ga=0$),\quad
$\Ahalf$ or $\Aminus$ ($\Ga=\pm\frac12$),\quad
$\Aroot$ ($\Ga=\pm1$),\quad
$\Aplus$ ($\Ga=\pm\frac32$),\quad
$\Art$ ($\Ga=\pm2$)}.
\]
Here, $\Aminus$ is used only for~$l_k$
and only if $\hbar_k\cdot l_k<0$; in all other cases, the
signs of~$l_k$ and $\hbar_k$ agree, so that we have $\hbar_k\cdot l_k\ge0$.

Recall that the kernel
\[*
\bN\bmod24\bA_2\subset\discr(24\bA_1)\cong(\Z/2)^{24}
\]
of the extension $\bN\supset24\bA_1$ is the \emph{extended binary Golay code}
$\CC_{24}$
(see~\cite{Conway.Sloane}).
The map $\supp\:\discr(24\bA_1)\to(\text{the power set of~$\Omega$})$
identifies codewords with
subsets of $\Omega$; then, $\CC_{24}$ is invariant
under complement and, in addition to~$\varnothing$ and~$\Omega$,
it consists of $759$ octads, $759$ complements thereof,
and $2576$ dodecads.

To simplify the notation, we identify the basis vectors $r_k$
(which are assumed fixed throughout) with their indices $k\in\Omega$. For a subset
$\CS\subset\Omega$, we let $\vv\CS:=\sum r$, $r\in\CS$,
and we abbreviate $\cw\CS:=\frac12\vv\CS\in\bN$
\emph{if $\CS\in\CC_{24}$ is a codeword}.

Now, as in the previous proofs, we treat the configurations
one by one.
We also give a more detailed description of each vector~$\hbar$ in terms of
the Golay code~$\Golay{24}$.

\subsection{Vector~\hlink{24A1}{1}}\label{s.24A1-1}
This set is replanted (see \autoref{s.replant}) to \config[Leech]{2}
in~$\Lambda$, see \autoref{s.Leech-2} below.
\gchecked

\def\porb#1{\obj{3}{\cluster_{#1}}}
\subsection{Vector~\hlink{24A1}{2}}\label{s.24A1-2}
We have $\ls|\stab\hh|=11520$ and
$\hh=\vv\CR$, where $\ls|\CR|=6$ and $\CR$ is a subset of an octad
$\CO:=\CR\cup\{r_1,r_2\}\in\Golay{24}$
(\cf. \autoref{s.24A1-11}).
Using the partition
\[*
\porb{c}:=\bigl\{\orb\in\oorb3\bigm|c\cap\supp\orb=\varnothing\bigr\},\quad
 c\in\pos_2(0),
\]
we show that $\Bii120[0,3]=\varnothing$.
There is but one discarded set, $\Lsub1$ of size $261$, that satisfies the
threshold~\eqref{eq.threshold};
this set can be described as
\[
\Lsub1=
\fF\cap\spn\bigl(\hh-4\cw{o},\cw{\Omega\sminus\CO},r_1,r_2\bigr)^\perp,
\sublabel1
\]
where
$\ls|o|=8$,
$\ls|o\cap\CO|=4$, $r_1\in o$, $r_2\notin o$.
(Till the end of this section, we reserve the notation~$o$ for a
\emph{codeword}, $o\in\Golay{24}$, intersecting the other fixed sets in a
certain prescribed way.
In each case,
one can easily check that the set of extra data used in the description
of an extremal geometric set is
unique up to $\stab\hh$.)
\gchecked

\def\porb#1{\obj{3}{\cluster_{#1}}}
\subsection{Vector~\hlink{24A1}{3}}\label{s.24A1-3}
We have $\ls|\stab\hh|=11520$ and $\hh=\cw\CO+\vv\CR$,
where $\CO$ is an octad and $\CR:=\{r_1,r_2\}\subset\CO$ is a $2$-element set.
Note that $\rank\fF=23$, as $\fF$ is annihilated by $\cw\CO-\vv\CR$; this
configuration can be re-embedded as \config{2}, acquiring a few extra conics
(increasing the rank to~$24$) and a few extra roots orthogonal to~$\hh$
(changing the combinatorial orbits).
Using the partition
\[*
\porb{c}:=\bigl\{\orb\in\oorb3\bigm|c\subset\supp\orb\bigr\},\quad
 c\in\pos_2\bigl(\tfrac12\bigr),
\]
we show that $\Bii126[0,3]=\varnothing$.
Three of the discarded sets satisfy~\eqref{eq.threshold}; they are
of the form
\[*
\fL=\fF\cap\spn\bigl(\hh-4r_1,\hh-4r_2,\cw{\Omega\sminus\CO},v\bigr)^\perp,
\]
where the fourth vector~$v$ is as follows:
\begin{alignat}3
\Lmax2:\quad&v=s,&\qquad&s\in\Omega\sminus\CO,
\maxlabel2\\
\Lmax3:\quad&v=\cw{o}-r_1,&&\mbox{$\ls|o|=8$, $\ls|o\cap\CO|=2$},&\quad&
\mbox{$r_1\in o$, $r_2\in o$},
\maxlabel3\\
\Lsub2:\quad&v=\cw{o}-r_1,&&\mbox{$\ls|o|=8$, $\ls|o\cap\CO|=4$},&&
\mbox{$r_1\in o$, $r_2\notin o$}.
\sublabel2
\end{alignat}
(As above, $o\in\Golay{24}$ and the subscript in the notation refers to
the cardinality.)
\gchecked

\def\porb#1{\obj{3}{\cluster_{#1}}}
\def\morb#1{\obj{3}{\mega_{#1}}}
\subsection{Vector~\hlink{24A1}{4}}\label{s.24A1-4}
We have $\ls|\stab\hh|=20160$ and $\hh=\cw\CO+r$, where $\CO$
is a codeword of length~$16$ and $r\in\CO$.
Using the (mega-)clusters
\[*
\porb{c}:=\bigl\{\orb\in\oorb3\bigm|c\subset\supp\orb\bigr\},\quad
\morb{k}:=\bigl\{\porb{c}\bigm|c\ni k\bigr\},\quad
 k\in\pos(0),\ c\in\pos_2(0),
\]
we show that
$\Bii124[3]=\varnothing$ and
$\defect[\fL;3]>210$ for each set $\fL\in\Bii85[0]$.
Two of the discarded sets satisfy the threshold~\eqref{eq.threshold}; they
are of the form
\[*
\fL=\fF\cap\spn\bigl(\hh-4r,\cw{\Omega\sminus\CO},s,v\bigr)^\perp,\quad
\]
where $s\in\Omega\sminus\CO$ and the fourth vector~$v$ is as follows:
\begin{alignat}3
\Lmax4:\quad&v=t,&\qquad&t\in\Omega\sminus\CO,\ t\ne s,
\maxlabel4\\
\Lsub3:\quad&v=\cw{o}-r,&&\mbox{$\ls|o|=\phantom{0}8$, $\ls|o\cap\CO|=4$},&\quad&
\mbox{$r\in o$, $s\notin o$},
\sublabel3\\
\Lmisc3:\quad&v=\cw{o}-r,&&\mbox{$\ls|o|=12$, $\ls|o\cap\CO|=6$},&&
\mbox{$r\notin o$, $s\notin o$}.
\misclabel3{249}
\end{alignat}
(The third set $\Lmisc3$ is discussed in \autoref{ex.249} below.)
\gchecked

\def\porb#1{\obj{2}{\cluster_{#1}}}
\def\pporb#1{\obj{3}{\cluster'_{#1}}}
\subsection{Vector~\hlink{24A1}{5}}\label{s.24A1-5}
We have $\ls|\stab\hh|=2304$ and $\hh=\cw\CO+\vv\CR$, where
$\CO$ is an octad and $\CR\subset\Omega\sminus\CO$ is a
$4$-element set such that there is an octad $o\supset\CR$ with the property
$\ls|o\cap\CO|=4$ (\cf. \autoref{s.24A1-12}).
Using the partitions
\begin{alignat*}2
\porb{c}&:=\bigl\{\orb\in\oorb2\bigm|c\subset\supp\orb\bigr\},\quad
 &c&\in\pos_2(2),\\
\pporb{k}&:=\bigl\{\orb\in\oorb3\bigm|k\notin\supp\orb\bigr\},\quad
 &k&\in\pos(2),
\end{alignat*}
we show
$\Bii53[2]=\Bii51[3]=\varnothing$ and
$\defect[\fL;3]>66$ for each set $\fL\in\Bii14[0]$.
\gchecked

\def\porb#1{\obj{2}{\cluster_{#1}}}
\def\pporb#1{\obj{4}{\cluster'_{#1}}}
\def\morb#1{\obj{2}{\mega_{#1}}}
\def\mmorb#1{\obj{4}{\mega'_{#1}}}
\subsection{Vector~\hlink{24A1}{6}}\label{s.24A1-6}
We have
$\ls|\stab\hh|=11520$ and $\hh=\cw\CO+\vv\CR$,
where
$\CO$ is a codeword of length~$16$ and
$\CR\subset\Omega\sminus\CO$ is a $2$-element set.
Using the (mega-)clusters
\begin{alignat*}4
\porb{k,n}&:=\bigl\{\orb\in\oorb2\bigm|k,n\in\supp\orb\bigr\},\quad
 &\morb{k}&:=\bigl\{\porb{k,n}\bigr\},\quad
 &k&\in\pos(0),\ &n&\in\pos(2),\\
\pporb{c}&:=\bigl\{\orb\in\oorb4\bigm|c\subset\supp\orb\bigr\},\quad
 &\mmorb{k}&:=\bigl\{\pporb{c}\bigm|c\ni k\bigr\},\quad
 &k&\in\pos(0),\ &c&\in\pos_2(0),
\end{alignat*}
we show that
$\Bii47[2]=\Bii49[4]=\varnothing$ and
$\defect[\fL;4]>130$ for each set $\fL\in\Bii80[0]$.
\gchecked

\def\porb#1{\obj{3}{\cluster_{#1}}}
\subsection{Vector~\hlink{24A1}{7}}\label{s.24A1-7}
We have $\ls|\stab\hh|=336$ and $\hh=\cw\CO+\vv\CR+r$, where $\CO\ni r$ is an
octad and $\CR\subset\Omega\sminus\CO$ is a $2$-element set.
Using the partition
\[*
\porb{k}:=\bigl\{\orb\in\oorb3\bigm|c\in\supp\orb\bigr\},\quad
 k\in\pos(2),
\]
we show that $\Bii15[3]=\varnothing$.
Besides, $\Bii5[4]=\Bii40[0,5]=\varnothing$.
\gchecked

\def\porb#1{\obj{2}{\cluster_{#1}}}
\def\pporb#1{\obj{3}{\cluster'_{#1}}}
\def\ppporb#1{\obj{4}{\cluster''_{#1}}}
\subsection{Vector~\hlink{24A1}{8}}\label{s.24A1-8}
We have $\ls|\stab\hh|=660$ and $\hh=\cw\CO+r+s$, where $\CO\ni r$ is a
dodecad and $s\in\Omega\sminus\CO$.
Using the partitions
\begin{alignat*}3
\porb{k}&:=\bigl\{\orb\in\oorb2\bigm|k\in\supp\orb\bigr\},\quad
 &k&\in\pos(0)&\quad&\mbox{($2$-fold)},\\
\pporb{k}&:=\bigl\{\orb\in\oorb3\bigm|k\in\supp\orb\bigr\},\quad
 &k&\in\pos(0)&&\mbox{($3$-fold)},\\
\ppporb{k}&:=\bigl\{\orb\in\oorb4\bigm|l_k\cdot\hh_k=-\tfrac12\bigr\},\quad
 &k&\in\pos\bigl(\tfrac12\bigr)&&\mbox{($1$-fold)},
\end{alignat*}
we show that
$\Bii38[2]=\Bii40[3]=\Bii11[0,4]=\varnothing$.
\gchecked

\def\porb#1{\obj{1}{\cluster_{#1}}}
\def\pporb#1{\obj{2}{\cluster'_{#1}}}
\def\mmorb#1{\obj{2}{\mega'_{#1}}}
\def\ppporb#1{\obj{4}{\cluster''_{#1}}}
\subsection{Vector~\hlink{24A1}{9}}\label{s.24A1-9}
We have $\ls|\stab\hh|=432$ and $\hh=\cw\CO+\vv\CR$, where $\CO$ is a dodecad
and $\CR\subset\Omega\sminus\CO$ is a $3$-element set.
Using the partitions
\begin{alignat*}3
\porb{k}&:=\bigl\{\orb\in\oorb1\bigm|k\in\supp\orb\bigr\},\quad
 &k&\in\pos(2)&\quad&\mbox{($1$-fold)},\\
\pporb{k,n}&:=\bigl\{\orb\in\oorb2\bigm|
 k\in\supp\orb,\ n\notin\supp\orb\bigr\},\quad
 &k&\in\pos(0),\ n\in\pos(2)
 &\quad&\mbox{($2$-fold)},
\end{alignat*}
and, grouping the latter into mega-clusters
$\mmorb{k}:=\{\pporb{k,n}\}$, $k\in\pos(0)$,
we show that
$\Bii13[1]=\Bii72[2]=\varnothing$ and
$\defect[\fL;2]>93$ for each set $\fL\in\Bii20[0]$.
In addition, we have
$\Bii3[4]=\Bii8[5]=\varnothing$.
\gchecked

\def\porb#1{\obj{2}{\cluster_{#1}}}
\subsection{Vector~\hlink{24A1}{11}}\label{s.24A1-11}\label{s.6D4-1}
We have $\ls|\stab\hh|=2160$ and
$\hh=\vv\CR$, where $\ls|\CR|=6$ and $\CR$ is \emph{not} a subset of an octad
(\cf. \autoref{s.24A1-2}).
Using the partition
\[*
\pporb{k}:=\bigl\{\orb\in\oorb2\bigm|k\notin\supp\orb\bigr\},\quad
 k\in\pos(2)\quad\text{($2$-fold)},
\]
we show that $\Bii24[ 2 ]=\varnothing$.
\gchecked

\def\porb#1{\obj{2}{\cluster_{#1}}}
\def\pporb#1{\obj{5}{\cluster_{#1}}}
\subsection{Vector~\hlink{24A1}{12}}\label{s.24A1-12}
We have $\ls|\stab\hh|=192$ and $\hh=\cw\CO+\vv\CR$, where
$\CO$ is an octad and $\CR\subset\Omega\sminus\CO$ is a
$4$-element set such that there is \emph{no} octad $o\supset\CR$ with the
property
$\ls|o\cap\CO|=4$ (\cf. \autoref{s.24A1-5}).
Using the partitions
\begin{alignat*}2
\porb{c}&:=\bigl\{\orb\in\oorb2\bigm|c\subset\supp\orb\bigr\},&\quad
 c&\in\pos_2(2),\\
\pporb{k}&:=\bigl\{\orb\in\oorb5\bigm|k\notin\supp\orb\bigr\},&\quad
 k&\in\pos(2),
\end{alignat*}
we show that $\Bii32[0,2]=\Bii31[5]=\varnothing$.
\gchecked

This case concludes the proof of \autoref{th.24A1}.
\endproof

\setlattice{Leech}
\section{The Leech lattice $\Lambda$\pdfstr{}{ \rm(\cf. \cite{DIO})}}\label{S.Leech}

The last lattice to be considered is the root-free \emph{Leech
lattice}~$\Lambda$.

\theorem\label{th.Leech}
With two exceptions \rom(up to automorphism\rom),
one has $\ls|\fL|\le260$ for any square~$12$ vector
$\hbar\in\Lambda$ and any geometric set $\fL\subset\fF(\Lambda,\hbar)$.
The exceptions are\rom:
\roster*
\item
one set $\Lmax5$ of size~$285$, see \eqref{eq.Lmax.5} below\rom;
\item
one set $\Lsub4$ of size~$261$, see \eqref{eq.Lsub.4} below.
\endroster
\endtheorem

\proof
According to Theorem 28 in \cite[Chapter~10]{Conway.Sloane},
any nonzero class $[\hbar]\in\Lambda/2\Lambda$ is represented
by a unique pair $\pm a\in\Lambda$,
{\4where $a^2\in\{4,6\}$ or $a^2=8$ and $a$ is an element of a fixed coordinate
frame.}
Since $a^2=\hbar^2\bmod4$ and $\Lambda$ is positive definite
and root free, for $\hbar^2=12$ we have
{\4(up to a basis change in $\Z a+\Z b$)}
either
\roster
\item\label{a=4}
$\anchor{1}\hbar = a + 2b$, where $a^2 = b^2 = 4$ and $a \cdot b = -2$
(\emph{type $6_{22}$} in \loccit.), or
\item\label{a=8}
$\anchor{2}\hbar = a + 2b$, where $a^2= 8$, $b^2 = 4$, and $a \cdot b = -3$
(\emph{type $6_{32}$} in \loccit.)
\endroster
Besides, a pair $a,b\in\Lambda$ as in item~\iref{a=4} or \iref{a=8}
is unique up to $\OG(\Lambda)$.
Thus, there are two $\OG(\Lambda)$-orbits of square~$12$ vectors
$\hbar\in\Lambda$ (see Theorem~29 in \loccit.)

\def\hmodtwo{\pdfstr{h = a mod 2L}{$\hh=a\bmod2\Lambda$}}

In each case, there is a unique orbit $\borb_1$ and, obviously, all
combinatorial orbits are single, $\Orb=\single$; therefore, we proceed as in
\autoref{s.single}, computing iterated index~$2$ subgroups and trying to find
large subsets $\fL\subset\fF$ of rank $\rank\fL\le20$.

\subsection{Vector \hlink{Leech}1: \hmodtwo, $a\sp2=4$}\label{s.Leech-1}
Let $F:=\spn_{\Z}\fF$. Then,
\[*
\ls|\OG_\hh(F)|=55180984320,\qquad
\ls|\fF|=891,\qquad
\rank F=23,\qquad
\discr F=\Cal{U}\oplus\bigl\<\tfrac74\bigr\>,
\]
and all index~$8$ extensions of $F\oplus\Z a$, $a^2=4$, are isomorphic and
have vector~$\hh$ of the same
type~\hlink{Leech}1. Note that $F$ is not primitive in~$\Lambda$ and
$\OG_\hh(\Lambda)$ induces an index~$6$ subgroup of $\OG_\hh(F)$. To
reduce the number of intermediate classes, we work with the smaller
lattice~$F$ and larger symmetry group $\OG_\hh(F)$; then, for each set
$\fL\subset\fF$ found, we analyze the lattice
$S:=\spn\fL\cap F$ itself as well as all
\emph{abstract} finite index extensions $\bar S\supset S$.

The computation results in four saturated sets $\fL\subset\fF$ of rank
$\rank\fL\le20$ and size $\ls|\fL|\ge261$: their sizes are $297$, $285$,
$279$, and $261$. None of the lattices $S$ as above has a nontrivial
\emph{root-free} finite
index extension $\bar S\supset S$ satisfying $\hh\in\bar S\dual$, and only
the two sets as in the statement of \autoref{th.Leech}
are geometric. These two sets can be described in terms of square~$4$ vectors
in~$\Lambda$.
To this end,
consider the lattice $U:=\Z a+\Z b+\Z c+\Z d$ with the Gram matrix
\[*
\bmatrix
  4&-2& 1& 1\\
 -2& 4&-2& 1\\
  1&-2& 4&-2\\
  1& 1&-2& 4\endbmatrix
\]
and let $V:=U+\Z v$, $v^2=4$, be its extension such that
$v\cdot a=v\cdot b=0$ and the
other two products are as follows:
\begin{alignat}3
&  v\cdot c=0,&&\quad v\cdot d=2
&&\quad\mbox{for $\Lmax5$},
\maxlabel5\\
&  v\cdot c=1,&&\quad v\cdot d=1
&&\quad\mbox{for $\Lsub4$}.
\sublabel4
\end{alignat}
Up to $\OG_\hh(\Lambda)$, there is a unique isometry $V\into\Lambda$ such
that
$a+2b\mapsto\hh$,
and the set in question is
\[*
\fF(\hh)\cap\bigl((\Z\hh\oplus V^\perp)\otimes\Q\bigr).
\]

\remark
The other two large sets can also be described in a similar manner. For the
largest one, we can take
$v\cdot a=v\cdot c=1$ and $v\cdot b=v\cdot d=-2$, whereas for the set of
size~$279$ the whole Gram matrix should be modified.
\endremark

\subsection{Vector \hlink{Leech}2: \hmodtwo, $a\sp2=8$}\label{s.Leech-2}
Letting, as above, $F:=\spn_\Z\fF$, we have
\[*
\OG_\hh(\Lambda)=\OG_\hh(F)=M_{24},\qquad
\ls|\fF|=759,\qquad
\rank F=24.
\]
Enumerating iterated index~$2$ subgroups, we find that each subset
$\fL\subset\fF$ of size $\ls|\fL|\ge261$ has rank $\rank\fL\ge21$ and,
therefore, is not geometric.
\endproof

\section{Proofs and examples}\label{S.proofs}

In this concluding section we complete the proof of \autoref{th.main} and
discuss a few interesting examples discovered in the course of the
computation.

\subsection{Proof of \autoref{th.main}}\label{proof.main}
By \autoref{th.existence} (and the Riemann--Roch theorem),
maximal (with respect to inclusion) deformation families of smooth sextic
surfaces $X\subset\Cp4$ whose generic member has
a prescribed dual adjacency graph~$\Gamma$ of conics are classified by
the oriented isomorphism classes of
triples $h\in{\NS}\subset\bL$ such that
\roster*
\item
$h^2=6$ and $\Fn_2(\NS,h)\cong\Gamma$;
\item
$h\in\NS\subset\bL$ satisfies
conditions~\iref{i.hyperbolic}--\iref{i.isotropic} in \autoref{s.homology},
and
\item
$\NS$ is rationally generated by~$h$ and $\Fn_2(\NS,h)$.
\endroster
(The last condition makes the N\'{e}ron--Severi lattice minimal and, hence, the
family maximal.) The construction of \autoref{s.S-lattice} reduces the study
of pairs ${\NS}\ni h$ \emph{admitting} a primitive embedding to~$\bL$ to the
study of geometric sets $\fL\cong\Gamma$ in $12$-polarized Niemeier lattices.
Finally, Theorems~\ref{th.<=8}, \ref{th.12A2}, \ref{th.24A1},
and~\ref{th.Leech} show that there are but nine geometric sets~$\fL$ of size
$\ls|\fL|\ge261$:
\roster*
\item
$\Lmax1$\,\eqref{eq.Lmax.1},
$\Lmax2$\,\eqref{eq.Lmax.2},
$\Lmax3$\,\eqref{eq.Lmax.3},
$\Lmax4$\,\eqref{eq.Lmax.4},
$\Lmax5$\,\eqref{eq.Lmax.5} of size~$285$ and
\item
$\Lsub1$\,\eqref{eq.Lsub.1},
$\Lsub2$\,\eqref{eq.Lsub.2},
$\Lsub3$\,\eqref{eq.Lsub.3},
$\Lsub4$\,\eqref{eq.Lsub.4} of size~$261$.
\endroster
Using the \texttt{GRAPE} package
\cite{GRAPE:nauty,GRAPE:paper,GRAPE} in \GAP~\cite{GAP4}, one can easily
establish that all five $\LMAX_{\sizemax}^*$ sets are pairwise isomorphic as
abstract graphs, and so are all four $\LSUB_{\sizesub}^*$ sets.
Furthermore, for each of the nine sets~$\fL$ one has $\spn_\Z\fL=\spn\fL$
and $\Hyp(\fL)=\varnothing$, see~\eqref{eq.Hyp}.
Thus, we obtain but two N\'{e}ron--Severi lattices,
$N_{285}:=\hyp0(\LMAX_{\sizemax}^*)$ and
$N_{261}:=\hyp0(\LSUB_{\sizesub}^*)$, which are both of rank~$20$
and type~$\I$.
In particular, it follows that the corresponding sextic surfaces are
projectively rigid and contain no lines, hence all conics are irreducible.

There remains to apply Nikulin's theory~\cite{Nikulin:forms} to classify the
primitive embeddings $h\in N_n\into\bL$, $n=285$ or $261$.
Using \texttt{GRAPE} again, we show that the canonical homomorphisms
\[*
\OG_h(N_n)=\Aut\Fn_2(N_n,h)\to\Aut\discr N_n
\]
are surjective and, moreover, any involution $\Ga\in\Aut\discr N_n$ lifts to
an involution $\tilde\Ga\in\OG_h(N_n)$. If $n=285$, the transcendental
lattice $T:={\NS}^\perp\subset\bL$ is unique in its genus, and we obtain a
single (formally, up to complex conjugation) surface
\begin{align}
\quartic1\:\quad&T=[ 6, 0, 20 ]\qlabel1{X_{285}}.
\end{align}
(We use the inline notation $[a,b,c]$ for a rank~$2$ lattice $\Z u+\Z v$,
$u^2=a$, $u\cdot v=b$, $v^2=c$.) If $n=261$, the genus of~$T$ has two
isomorphism classes, giving rise to two sextics
(not isomorphic even as abstract $K3$-surfaces)
\begin{align}
\quartic2\:\quad&T=[ 8, 0, 18 ]\qlabel2{X_{261}},\\
\quartic3\:\quad&T=[ 2, 0, 72 ]\qlabel3{Y_{261}}
\end{align}
sharing the same combinatorial configuration of conics~$\LSUB_{\sizesub}^*$.

Since each of the transcendental lattices~$T$ has an involutive
orientation reversing
automorphism $\Gb$, which, as explained above, can be matched with an
involutive automorphism $\tilde\Ga\in\OG_h(N_n)$, all three surfaces are
real. The last transcendental lattice has a root~$r$;
hence, $\quartic3$ (and only this surface) has a real structure (\viz.
the one that
induces the reflection $-\tr_r$) with respect to which all conics are real
(\cf. \cite[Lemma~3.8]{DIS}). One can easily check that, for each conic
$c\in\Fn_2(N_{261})$, there is another conic~$c'$ such that $c\cdot c'=1$;
hence, each conic has a real point.
\qed

\subsection{Other examples}\label{s.examples}
Altogether, we have found $71$ combinatorial configurations of more than
$200$ conics. The ten largest
configurations~$\graph$
are presented in
\autoref{tab.large},
where we show the size $\ls|\graph|$ of~$\graph$,
the size $\ls|\Aut\graph|$
of its automorphism group, the transcendental lattice(s)~$T$ \emph{provided
that all conics are irreducible}, and a reference to the description of
reducible conics, if any.
In addition to \autoref{tab.large},
there are $61$ configurations representing all \emph{odd} counts
between $201$ and~$225$. We do not assert that any of these lists is
complete.
\table
\caption{Known large configurations~$\graph$ of conics}\label{tab.large}
\hbox to\hsize{\hss\vbox{%
\halign{\quad\strut\hss$#$\hss&\quad\hss$#$\hss&\quad$#$\hss&\quad#\hss\quad\cr
\noalign{\hrule\vspace{2pt}}
\ls|\graph|&\ls|\Aut\graph|&T:=\NS(X)^\perp&Reducible?\cr
\noalign{\vspace{1pt}\hrule\vspace{2pt}}
285&2880&[ 6, 0, 20 ]\cr
261& 288&[ 2, 0, 72 ], [ 8, 0, 18 ]\cr
249& 144&[ 2, 0, 78 ], [ 6, 0, 26 ]&\autoref{ex.249}\cr
243& 144&[ 12, 0, 14 ]\cr
237&2880&[ 10, 0, 20 ]\cr
237&  96&[ 10, 4, 20 ]\cr
237&  48&[ 2, 0, 88 ], [ 10, 2, 18 ]\cr
231& 288&&\autoref{ex.231}\cr
231&  72&[ 8, 4, 26 ]\cr
229&  43&[ 10, 4, 20 ]\cr
\noalign{\vspace{1pt}\hrule}
\crcr}}\hss}
\endtable

It is worth mentioning that, unlike large configurations of lines
(see~\cite{degt:lines}), large configurations of conics are irregular: the
group $\Aut\graph$ tends to have many orbits.

Most large configurations appear in sextic surfaces of type~$\I$ only,
\ie, there are no lines (or other curves of odd projective degree)
and all conics are irreducible. The few known large configurations that
\emph{may} contain reducible conics are discussed in the examples below;
remarkably, all $\det$-extremal sextics found in~\cite{degt:singular.K3}
appear on the list.
We
indicate the
the numbers of conics in the form
\[*
\text{(total)}=\text{(irreducible)}+\text{(reducible)}.
\]

\example\label{ex.249}
The third largest configuration
of conics that we have observed is the triple of
pairwise isomorphic geometric sets
\roster*
\item
$\Lmisc1$\,\eqref{eq.Lmisc.1},
$\Lmisc2$\,\eqref{eq.Lmisc.2},
$\Lmisc3$\,\eqref{eq.Lmisc.3} of size~$249$.
\endroster
This time, $\Hyp(\LMISC_{249}^*)=\{\kappa\}$ is a one-element set, and
\autoref{prop.embedding} asserts that both $\hyp0(\LMISC_{249}^*)$ and
$\hyp\kappa(\LMISC_{249}^*)$ admit a primitive embedding to~$\bL$,
\ie, the same abstract configuration is
realized both by all irreducible conics and by (partially) reducible ones.
Arguing as in \autoref{proof.main}, we arrive at four surfaces
(where the respective transcendental lattice $T$ is incorporated in the notation):
\roster*
\item
$X_{249}([ 6, 0, 26 ])$, $Y_{249}([ 2, 0, 78 ])$: $249$ irreducible
conics, and
\item
$X_{42,60}([ 6, 3, 8 ])$, $Y_{42,60}([ 2, 1, 20 ])$:
$42$ lines and
$249=60+189$
conics.
\endroster
All four are real, and the two $Y_*$-surfaces admit real structures with
respect to which all lines (if any) and conics are real. Note also that
$X_{42,60}$ and $Y_{42,60}$ are the two discriminant minimizing singular
(in the sense of the maximal Picard rank)
smooth sextics discovered in~\cite{degt:singular.K3} ($6_{42}$ in the
notation thereof),
and $42$ is the maximal number of lines in a smooth sextic surface
(see~\cite{degt:lines}).
\endexample

\remark\label{rem.189}
One can easily show that the maximal number of \emph{reducible} conics in a
smooth sextic surface~$X$ is~$189$, and this number is only realized in the
irreducible $1$-parameter family $\Psi_{42}$ of sextics
with the maximal number~$42$ of
lines (see~\cite{degt:lines}).

Indeed, according to Proposition~2.12 in \loccit. and a remark thereafter,
the maximal valency of a vertex $v\in\Fn_1X$ is $\val v\le9$, {\em unless $v$
is part of a triangle \rom(cycle of length~$3$\rom) and $\val v=11$.}
By Lemma~8.3, any two triangles in $\Fn_1X$ are disjoint, and by Lemma~8.4,
the total valency of the three vertices in a triangle is at most
$21+6=27$. Hence, the \emph{average} valency of a vertex is still at most~$9$, and
the number of reducible conics is at most $\frac12\cdot42\cdot9=189$.
\endremark

\example\label{ex.231}
Each member of the family~$\Psi_{42}$ (see \autoref{ex.249}) has
the same configuration of $42$ lines
and, hence, $189$ reducible conics.
In addition, a generic member of~$\Psi_{42}$ has $24$ irreducible conics,
and this combinatorial configuration can alternatively be realized by
$213$ irreducible conics. Apart from \autoref{ex.249}, there is at least one
more singular member with more conics:
\roster*
\item
$X_{42,42}([ 6, 0, 8 ])$ ($6'_{42}$ in~\cite{degt:singular.K3}):
$42$ lines and
$231=42+189$
conics.
\endroster
This configuration can not be realized by $231$ irreducible conics.
However, it
\emph{does} admit an alternative realization, with a different set of
reducible conics,
by another embedding of the same singular $K3$-surface
$X([ 6, 0, 8 ])$ to $\Cp4$:
\roster*
\item
$X_{36,87}([ 6, 0, 8 ])$ ($6'_{36}$ in~\cite{degt:singular.K3}):
$36$ lines and
$231=87+144$
conics.
\endroster
The third smooth embedding of the same surface is
\roster*
\item
$X_{38,70}([ 6, 0, 8 ])$ ($6_{38}$ in~\cite{degt:singular.K3}):
$38$ lines and
$223=70+153$
conics,
\endroster
and the remaining surface found in~\cite{degt:singular.K3} is
\roster*
\item
$X_{36,81}([ 8, 4, 8 ])$ ($6_{36}''$ in~\cite{degt:singular.K3}):
$36$ lines and
$225=81+144$
conics.
\endroster
\endexample

{
\let\.\DOTaccent
\def\cprime{$'$}
\bibliographystyle{amsplain}
\bibliography{degt}
}

\end{document}